\documentclass[11pt]{article}
\usepackage[a4paper, margin=2cm]{geometry}
\usepackage[T1]{fontenc}
\usepackage{amsmath, amssymb}
\usepackage{hyperref}
\usepackage{authblk}
\usepackage{booktabs}
\usepackage{multicol}
\usepackage{enumitem}
\usepackage{subcaption}
\usepackage{tikz}
\usepackage{tikz-3dplot}
\usetikzlibrary{3d,shapes,arrows,positioning,calc,patterns}
\usepackage{pgfplots}\pgfplotsset{compat=newest}
\usepgfplotslibrary{fillbetween,colormaps}
\usepackage{gensymb}
\usepackage{cleveref}
\crefname{figure}{Figure}{Figures}
\crefname{section}{Section}{Sections}
\usepackage[style=numeric, giveninits,doi=false,isbn=false,url=false,eprint=false]{biblatex}
\pgfplotsset{
colormap/plasma/.style={%
    /pgfplots/colormap={plasma}{%
      rgb=(0.050383, 0.029803, 0.527975)
      rgb=(0.186213, 0.018803, 0.587228)
      rgb=(0.287076, 0.010855, 0.627295)
      rgb=(0.381047, 0.001814, 0.653068)
      rgb=(0.471457, 0.005678, 0.659897)
      rgb=(0.557243, 0.047331, 0.643443)
      rgb=(0.636008, 0.112092, 0.605205)
      rgb=(0.706178, 0.178437, 0.553657)
      rgb=(0.768090, 0.244817, 0.498465)
      rgb=(0.823132, 0.311261, 0.444806)
      rgb=(0.872303, 0.378774, 0.393355)
      rgb=(0.915471, 0.448807, 0.342890)
      rgb=(0.951344, 0.522850, 0.292275)
      rgb=(0.977856, 0.602051, 0.241387)
      rgb=(0.992541, 0.687030, 0.192170)
      rgb=(0.992505, 0.777967, 0.152855)
      rgb=(0.974443, 0.874622, 0.144061)
      rgb=(0.940015, 0.975158, 0.131326)
    },
  },
}

\providecommand\bx{\boldsymbol{x}}
\providecommand\bX{\boldsymbol{X}}
\providecommand\bn{\boldsymbol{n}}
\providecommand\by{\boldsymbol{y}}
\providecommand\bu{\boldsymbol{u}}

\providecommand\bq{\boldsymbol{q}}
\providecommand\btau{\boldsymbol{\tau}}
\providecommand\baq{\boldsymbol{\acute{q}}}
\providecommand\bhv{\boldsymbol{\hat{v}}}
\providecommand\bhw{\boldsymbol{\hat{w}}}
\providecommand\bhq{\boldsymbol{\hat{q}}}
\providecommand\bbq{\boldsymbol{\breve{q}}}
\providecommand\bhbq{\boldsymbol{\hat{\breve{q}}}}
\providecommand\bdq{\boldsymbol{\delta{}q}}
\providecommand\bdhq{\boldsymbol{\delta{}\hat{q}}}
\providecommand\beff{\boldsymbol{f}}
\providecommand\bqd{\boldsymbol{\dot{q}}}

\providecommand\bnabla{\boldsymbol{\nabla}}
\providecommand\bcdot{\boldsymbol{\cdot}}
\providecommand\calF{\mathcal{F}}

\providecommand\calM{\mathcal{M}}
\providecommand\calR{\mathcal{R}}
\providecommand\calD{\mathcal{D}}

\providecommand\calI{\mathcal{I}}
\providecommand\balpha{\boldsymbol{\alpha}}
\providecommand\bzero{\boldsymbol{0}}

\definecolor{paperred}{RGB}{238,32,77}
\definecolor{paperorange}{RGB}{255,117,56}
\definecolor{paperyellow}{RGB}{255,170,29}
\definecolor{papergreen}{RGB}{28,172,120}
\definecolor{paperblue}{RGB}{31,117,254}
\definecolor{paperviolet}{RGB}{150,61,127}

\title{\texttt{ff-bifbox}: A scalable, open-source toolbox for bifurcation analysis of nonlinear PDEs}
\author[1]{Christopher M. Douglas}
\author[2]{Pierre Jolivet}
\affil[1]{Duke University, Durham, NC, USA, 27708}
\affil[2]{Sorbonne Universit\'{e}, CNRS, LIP6, Paris, France, 75252}
\begin{document}
\maketitle
%% ABSTRACT
\begin{abstract}
Nonlinear PDEs give rise to complex dynamics that are often difficult to analyze in state space due to their relatively large numbers of degrees of freedom, ill-conditioned operators, and changing spatial and parameter resolution requirements. This work introduces \texttt{ff-bifbox}: a new open-source toolbox for performing numerical branch tracing, stability/bifurcation analysis, resolvent analysis, and time integration of large, time-dependent nonlinear PDEs discretized on adaptively refined meshes in two and three spatial dimensions. Spatial discretization is handled using finite elements in FreeFEM, with the discretized operators manipulated in a distributed framework via PETSc. Following a summary of the underlying theory and numerics, results from three examples are presented to validate the implementation and demonstrate its capabilities. The considered examples, which are provided with runnable \texttt{ff-bifbox} code, include: a 3-D Brusselator system, a 3-D plate buckling system, and a 2-D compressible Navier--Stokes system. In addition to reproducing results from prior studies, novel results are presented for each system.
\end{abstract}

%%%%%%%%%%%%%%%%%%%%%%%%%%%
\section{Introduction}\label{sec:introduction}
%%%%%%%%%%%%%%%%%%%%%%%%%%%

Dynamical systems analysis is a powerful framework for gaining fundamental insight into the behavior of complex nonlinear systems. Through the systematic exploration of a model's phase space and the identification of its attractors, bifurcations, and other critical features, dynamical systems analysis has contributed deep phenomenological insights in many scientific fields and enabled revolutionary control strategies across diverse branches of engineering~\cite{strogatz_2018}. Much of this progress can be linked to the proliferation of open-source software like AUTO~\cite{doedel_1981,doedel_etal_2007}, MATCONT~\cite{dhooge_etal_2008}, PyDSTool~\cite{clewley_etal_2012}, and CoCo~\cite{dankowicz_schilder_2013} for numerical branch tracing and bifurcation analysis of generic, nonlinear ordinary differential equations (ODEs) and systems of nonlinear ODEs. However, the extension of such packages, which are designed for dense and relatively compact equation systems, to sparse spatiotemporal systems described by nonlinear partial differential equations (PDEs) remains a substantial challenge. Indeed, discretized systems of PDEs posed on non-trivial domains often yield very large ODE systems with $\mathcal{O}(10^6)$ degrees of freedom or more, requiring specialized algorithms beyond the generalist tools already available for standard ODEs.

In response to this challenge, a handful of open-source packages focused on PDE applications have been developed. For example, the pde2path~\cite{uecker_2021} and BifurcationKit.jl~\cite{veltz_2020} projects both offer comprehensive branch tracing and bifurcation analysis capabilities for smaller PDE systems on fixed meshes. Another notable package is LOCA~\cite{salinger_etal_2002}, now integrated with the Trilinos project~\cite{trilinos-website}, which provides generic branch tracing and bifurcation analysis capabilities for large-scale, parameter-dependent PDEs. There are also some excellent application-focused projects like oomph-lib~\cite{heil_hazel_2006}, StabFEM~\cite{fabre_etal_2019}, and BROADCAST~\cite{poulain_etal_2023}, which implement practical tools for many aspects of PDE analysis but do not focus on a unified, scalable framework for parametric bifurcation analysis of large PDE systems.

This article introduces \texttt{ff-bifbox}, a package of cross-compatible FreeFEM scripts designed explicitly for numerical continuation, stability/bifurcation analysis, resolvent analysis, and time integration of large, nonlinear PDEs discretized on adaptively-refined meshes in two and three spatial dimensions. The open-source project is built on top of FreeFEM~\cite{hecht_2012}, a finite-element software and scripting language with native two-dimensional adaptive meshing capabilities, and PETSc~\cite{balay_etal_2026}, a scalable scientific computing library with versatility for a huge range of different solvers and preconditioners. \texttt{ff-bifbox} implements robust, adaptive continuation routines leveraging efficient block factorizations that promote sparsity to reduce communication and memory requirements, and preserve the underlying structure of the discrete operators such that a preconditioner suitable for a given system can be efficiently reused throughout the library. Furthermore, solutions from the various \texttt{ff-bifbox} routines can all be directly exported for easy post-processing and visualization in ParaView~\cite{ahrens_etal_2005}.

The remainder of this article is organized as follows. \Cref{sec:theory} introduces essential theoretical aspects of bifurcation analysis and presents key algorithms as well as their implementation and structure in the \texttt{ff-bifbox} project. \Cref{sec:validation} provides example computations of continuum mechanics problems that demonstrate the use and validate the results of \texttt{ff-bifbox}. Finally, \cref{sec:conclusion} concludes the article with a summary of the current capabilities and limitations of \texttt{ff-bifbox} and describes some objectives for the project's future development. 

%%%%%%%%%%%%%%%%%%%%%%%%%%%
\section{Theory and Implementation} \label{sec:theory}
%%%%%%%%%%%%%%%%%%%%%%%%%%%

\texttt{ff-bifbox} is designed to handle autonomous\footnote{In its current iteration, \texttt{ff-bifbox} does support some capabilities for non-autonomous systems analysis that will not be discussed here.}, continuous-time dynamical systems defined by a nonlinear function $\calF$, which may be written in the form:
\begin{align}\label{eq:nonlinearPDE}
    \calF\left(\bq,\bqd;\balpha\right)=\calM\left(\bq;\balpha\right)\bqd+\calR\left(\bq;\balpha\right)=0,
\end{align}
where $\calM$ is the square mass operator and $\calR$ is the steady residual vector. Here, the state vector $\bq\left(\bx,t\right)$ is a function of space $\bx$ and time $t$, $\bqd$ is the time derivative of $\bq$, and $\balpha$ is a vector of parameters. For a more comprehensive introduction to dynamical systems analysis, as well as a more pedagogical exposition of the methods outlined below, the interested reader is referred to any of the excellent monographs by Govaerts~\cite{govaerts_2000}, Kuznetsov~\cite{kuznetsov_2023}, and Allgower \& Georg~\cite{allgower_georg_2003}, among others.

The \texttt{ff-bifbox} implementation requires the user to define the computational domain using a mesh file, a variational formulation for the system of interest using the FreeFEM scripting language, and a set of state variables, parameters, and, where needed, user-defined functions or macros for more specialized tasks. From this basic starting point, \texttt{ff-bifbox} supplies a collection of command-line accessible algorithms, file management macros, and post-processing routines supporting various aspects of dynamical systems analysis for large-scale PDEs. It should be remarked that the weak formulation provides considerable flexibility, supporting exotic boundary conditions (e.g., \cite{douglas_2024}), non-smooth or discontinuous systems (e.g., via discontinuous elements), and non-standard couplings and constraints, including multiphysics and coupled ODE–PDE models.

In the following, the Jacobian operator $\partial_{\bq}\calF|_{\bq,\balpha}\left(\cdot\right)$ is denoted as the partial derivative of $\calF$ with respect to $\bq$, i.e. $\frac{\partial\calF}{\partial\bq}\left(\bq;\balpha\right)$. The parameter gradient operator $\frac{\partial\calF}{\partial\balpha}\left(\bq;\balpha\right)=\partial_{\balpha}\calF|_{\bq,\balpha}\left(\cdot\right)$ is also defined. Similar abuse of notation is used to define the second parameter gradient operator $\partial_{\balpha\balpha}\calF|_{\bq,\balpha}\left(\cdot,\cdot\right)$, the Hessian operator $\partial_{\bq\bq}\calF|_{\bq,\balpha}\left(\cdot,\cdot\right)=\frac{\partial^2\calF}{\partial\bq\partial\bq}\left(\bq;\balpha\right)$, and others.

%%%%%%%%%%%%%%%%%%%%%%%%%%%
\subsection{Steady solutions}\label{ssec:fixedpoints}

Steady solutions to \cref{eq:nonlinearPDE}, also termed fixed points or equilibria, are denoted by $\bqd=\bzero$ and satisfy the definition,
\begin{align}\label{eq:equilibrium}
    \calF\left(\bq,\bzero;\balpha\right)=\calR\left(\bq;\balpha\right)=0.
\end{align}
With an initial condition and $\balpha$ given, solutions to \cref{eq:equilibrium} may be found using the \texttt{SNES} library from PETSc~\cite{balay_etal_2026}, which implements a variety of Newton-like solvers based on an initial guess $\bq_0$ and a sequence of iterations consisting of the steps,
\begin{subequations}\label{eq:Newton}
\begin{align}\label{eq:Newtona}
    \partial_{\bq}\calR|_{\bq_k,\balpha}\left(\bdq_{k}\right)&=\calR\left(\bq_k;\balpha\right),\\\label{eq:Newtonb}
    \bq_{k+1}&=\bq_k-\beff\left(\bdq_k\right),\qquad{}k=0,1,\hdots,
\end{align}
\end{subequations}
where $\partial_{\bq}\calR{}|_{\bq_k,\balpha}\left(\cdot\right)$ is the steady Jacobian operator evaluated at the point $(\bq_k,\balpha)$, $\bdq_k$ is a correction, $k$ is the index of the present iteration, and $\beff$ is some function used by \texttt{SNES} to control the update. This iterative sequence is repeated until a specified convergence or divergence threshold is met.

Once a single steady solution is identified, it is often desirable to extend that solution over some region of the parameter space along its equilibrium manifold or to identify other distinct solutions at the same parameter values. These aims, respectively termed ``branch tracing'' and ``branch switching'' can be accomplished using predictor--corrector methods and deflation methods.

\subsubsection{Predictor--corrector method for branch tracing}\label{sssec:tracing}

A robust approach to continuation of solutions along parameters involves treating both the state vector $\bq$ and an arbitrary parameter $\alpha\in\balpha$ as unknowns. This yields a predictor--corrector scheme that can pass over singularities of the Jacobian operator (where the convergence radius of Newton's method goes to zero)~\cite{keller_1978}. Here, the predictor is defined by a tangent step whose direction is given by the augmented null vector $\by\left(\bq,\balpha\right)=\left[\by_{\bq}\left(\bq,\balpha\right),y_\alpha\left(\bq,\balpha\right)\right]^T\in\ker\left(\left[\partial_{\bq}\calR|_{\bq,\balpha},\ \partial_{\alpha}\calR|_{\bq,\balpha}\right]^T\right)$. Provided $\partial_{\bq}\calR|_{\bq,\balpha}$ is non-singular, the null vector may be identified by setting $y_{\alpha}=-1$ to fix its orientation with respect to the continuation parameter $\alpha$ without loss of generality, and then solving the linear system,
\begin{align}\label{eq:MPtangent}
    \partial_{\bq}\calR|_{\bq,\balpha}\left(\by_{\bq}\right)&=\partial_{\alpha}\calR|_{\bq,\balpha},
\end{align}
to determine $\by_{\bq}$. The tangent predictor step is then taken to be $\left[\bq,\alpha\right]^T=\left[\bq,\alpha\right]^T+\frac{h}{\left\|\by\right\|_2}\left[\by_{\bq},-1\right]^T$, where $h$ is a numerical parameter controlling the orientation and magnitude of the predictor step and $\|\by\|_2=\sqrt{\|\by_{\bq}\|^2+1}$ is the standard 2-norm of $\by$. To maintain a consistent direction and enable robust branch tracing without the computational burden of excessively small steps, the adaptive steplength scaling procedure formulated by Allgower \& Georg~\cite[Ch. 6]{allgower_georg_2003} is used to automatically adjust $h$ at each step.

Following each predictor step, a sequence of constrained Newton-like corrector iterations are required to converge from the prediction to an acceptable solution without stepping back along the predictor. The implementation in \texttt{ff-bifbox} leverages a Moore--Penrose approach where each iteration involves the augmented linear system,
\begin{subequations}\label{eq:MPcorrector}
    \begin{align}\label{eq:MPcorrectora}
    \begin{bmatrix}
        \partial_{\bq}\calR|_{\bq_k,\balpha_k}\left(\cdot\right)&\partial_{\alpha}\calR|_{\bq_k,\balpha_k}\\
        \langle\by_{\bq}\left(\bq_k,\balpha_k\right),\cdot\rangle&-1
    \end{bmatrix}
    \begin{bmatrix}
        \bdq_k\\
        \delta\alpha_k
    \end{bmatrix}
    &=
    \begin{bmatrix}
        \calR\left(\bq_k;\balpha_k\right)\\0
    \end{bmatrix},\\
    \begin{bmatrix}
        \bq_{k+1}\\
        \alpha_{k+1}
    \end{bmatrix}&=\begin{bmatrix}
        \bq_k\\
        \alpha_k
    \end{bmatrix}-\beff\left(\begin{bmatrix}
        \bdq_k\\
        \delta\alpha_k
    \end{bmatrix}\right),\qquad{}k=0,1,\hdots
    \end{align}
\end{subequations}
This iterative sequence for the corrector is repeated until the specified convergence threshold is met, then the null vector $\by$ at the converged point is used as the predictor for the following step and the sequence is repeated.

In \texttt{ff-bifbox}, the corrector sequence \cref{eq:MPcorrector} utilizes the \texttt{SNES} library, leveraging an exact block Schur decomposition from \texttt{PCFIELDSPLIT} to efficiently resolve \cref{eq:MPcorrectora}. This block Schur decomposition is preferred over a monolithic solution of the augmented operator for two important reasons. First, the block Schur decomposition respects the sparsity and structure of $\partial_{\bq}\calR$, allowing preconditioners developed for the Jacobian (or approximate Jacobian) operator to be used instead of requiring specialized new preconditioners for the augmented system. Second, even if a monolithic approach were used, an additional calculation would be needed to update $\by_k$ after each corrector iteration. On the other hand, the exact block Schur decomposition involves intermediate solutions of \cref{eq:Newtona} and \cref{eq:MPtangent}, meaning that the tangent vector $\by_k$ is obtained automatically as a byproduct of each iteration (see~\cite{douglas_etal_2021b} for further details).

\subsubsection{Deflation technique for branch switching}\label{sssec:switching}

Deflation is a strategy to identify multiple distinct solutions of a nonlinear system at fixed parameter values. Formally, deflation techniques seek to modify the governing system so that known solutions are removed from the solution set, thereby potentially enabling convergence to additional roots that would otherwise remain hidden~\cite{brown_gearhart_1971}. Following Farrell \textit{et al.}~\cite{farrell_etal_2015}, one can define the deflation operator,
\begin{align}\label{eq:deflation}
    \calD_{p,a}\left(\bq,\bq_j\right)=\frac{\calI}{\left\|\bq-\bq_j\right\|^p}+a\calI,
\end{align}
where $\calI{}$ is the identity operator, $\bq_j$ is a known solution, $p$ is the deflation order parameter, and $a$ is a shift parameter. Therefore, modifying \cref{eq:Newtona} through application of $\calD_{p,a}$ to obtain the deflated residual $\calD{}_{p,a}\left(\bq,\bq_0\right)\calR\left(\bq;\balpha\right)=0$ formally prevents $\bq_0$ from being a valid root of the deflated system, provided $p$ is large enough.

In contrast to the original approach of~\cite{farrell_etal_2015}, the implementation in \texttt{ff-bifbox} relies on an algebraic simplification introduced by Adler \textit{et al.}~\cite{adler_etal_2017}. Their key observation was that each Newton step with the deflated system differs from the original step by only a scalar factor. More specifically, with $\calD_{p,a}$ given by \cref{eq:deflation}, the deflated correction step is equivalent to \cref{eq:Newtonb}, but with $\beff\left(\bdq_k\right)$ replaced by,
\begin{align}
    \beff\left(\prod_j\left[\frac{1+a\left\|\bq_k-\bq_j\right\|^p}{1+a\left\|\bq_k-\bq_j\right\|^p-p\left\|\bq_k-\bq_j\right\|^{-2}\left\langle\bq_k-\bq_j,\bdq_k\right\rangle}\right]\bdq_k\right),
\end{align}
where the product operator serves to allow deflation of multiple solutions.

This formulation avoids any construction and linearization of $\calD_{p,a}$ while preserving its effect on the nonlinear iterations of the deflated system though a simple modification of the step magnitude and sign. Thus, the deflation approach adopted in \texttt{ff-bifbox} only involves evaluating a few additional inner products between arrays that are already computed in the original system \cref{eq:Newton}. This is efficiently implemented using \texttt{PCSHELL} without requiring additional linear solves or constraints.

%%%%%%%%%%%%%%%%%%%%%%%%%%%
\subsection{Stability of fixed-point solutions}\label{ssec:stability}

The stability of a steady solution to \cref{eq:equilibrium} is determined by the time evolution of superimposed perturbations. Using the Laplace transform, any such perturbation $\baq$ may generically be expressed as a superposition of spatial modes $\bhq$ of the form,
\begin{align}\label{eq:superposition}
    \baq\left(\bx,t\right)=\sum\left[\bhq\left(\bx\right)\exp\left(\lambda{}t\right)+\bhq^*\left(\bx\right)\exp\left(\lambda{}^*t\right)\right],
\end{align}
where $\left(\cdot\right)^*$ denotes complex conjugation and $\lambda$ is a complex scalar representing the Laplace variable. In general, the modes $\bhq$ are neither independent nor orthogonal. However, by restricting the stability analysis to infinitesimally small disturbances and considering their behavior at only asymptotically large times, the modes in \cref{eq:superposition} become both independent and orthogonal to each other. That is, they become equivalent to the generalized eigenmodes of \cref{eq:nonlinearPDE} when linearized about a solution to \cref{eq:equilibrium}. Therefore, the time-asymptotic stability of a fixed point to infinitesimal perturbations is determined by the generalized eigenvalue spectrum of the linearized operator evaluated at the fixed point,
\begin{align}\label{eq:stability}
    \left\langle\bhq^{\dagger{}},\left[\lambda{}\calM\left(\bq;\balpha\right)+\partial_{\bq}\calR|_{\bq,\balpha}\left(\cdot\right)\right]\bhq\right\rangle&=0.
\end{align}
Here, $\langle\cdot,\cdot\rangle$ denotes the standard Hermitian inner product, $\lambda{}\calM\left(\bq;\balpha\right)+\partial_{\bq}\calR|_{\bq,\balpha}\left(\cdot\right)$ is the linearized operator, $\bhq$ and $\bhq^{\dagger}$ represent, respectively, a direct (right) and adjoint (left) eigenmode pair, and $\lambda=\sigma+\mathrm{i}\omega$ is the corresponding eigenvalue with growth rate $\Re\left\{\lambda\right\}=\sigma$ and frequency $\Im\left\{\lambda\right\}=\omega$. If all eigenvalues satisfy $\sigma<0$, the system is defined to be asymptotically stable; conversely, the system is unstable if any eigenvalue has $\sigma>0$. This notably implies the existence of a critical state when $\sigma=0$. Such states correspond to bifurcations, which will be given further attention in \cref{ssec:bifurcations}.

Solutions to \cref{eq:stability} are identified numerically using the \texttt{EPS} library in SLEPc~\cite{hernandez_etal_2005}, leveraging a Krylov--Schur method based on the shift-and-invert technique. Importantly, this approach requires a choice of shift(s) for efficient resolution of the problem. Since the typical eigenvalues of interest in bifurcation problems fall near the line of neutral stability ($\sigma\sim0$), a natural choice is to analyze multiple shifts distributed along the imaginary axis to cover all relevant timescales of the considered dynamical system.

%%%%%%%%%%%%%%%%%%%%%%%%%%%
\subsection{Bifurcations of fixed-point solutions}\label{ssec:bifurcations}

Generic local bifurcations of fixed points with codimension 1 satisfy \cref{eq:equilibrium} under a criticality condition that $\sigma=0$ in \cref{eq:stability}. This criticality requirement can be enforced using several distinct approaches, as discussed by Govaerts~\cite{govaerts_2000}. In \texttt{ff-bifbox}, the constraint is enforced similarly for saddle--node (i.e. fold) bifurcations, branch points (e.g. pitchfork bifurcations), and Hopf bifurcations using a minimally augmented formulation with a scalar criticality variable. In the most general case of a Hopf bifurcation, the scalar criticality variable $g$ may be formally defined as,
\begin{align}\label{eq:criticalityvariable}
g=\left\langle\bhq^{\dagger{}},\left[\mathrm{i}\omega\calM\left(\bq;\balpha\right)+\partial_{\bq}\calR|_{\bq,\balpha}\left(\cdot\right)\right]\bhq\right\rangle,
\end{align}
which vanishes when the direct and adjoint eigenmodes satisfy \cref{eq:stability} with $\sigma=0$. However, \cref{eq:criticalityvariable} is only a formal definition of $g$. In practice, the criticality variable is instead determined using the bordered linear systems,
\begin{subequations}\label{eq:minaug}
\begin{align}
    \begin{bmatrix}
        -\left[\mathrm{i}\omega\calM\left(\bq;\balpha\right)+\partial_{\bq}\calR|_{\bq,\balpha}\left(\cdot\right)\right] & \calM\left(\bq;\balpha\right)\bhw\\
        \langle\calM\left(\bq;\balpha\right)\bhv,\cdot\rangle&0
    \end{bmatrix}
    \begin{bmatrix}
        \bhq\\g
    \end{bmatrix}
    &=
    \begin{bmatrix}
        0\\1
    \end{bmatrix},\\
    \begin{bmatrix}
        -\left[\mathrm{i}\omega\calM\left(\bq;\balpha\right)+\partial_{\bq}\calR|_{\bq,\balpha}\left(\cdot\right)\right]^{\dagger} & \calM\left(\bq;\balpha\right)\bhv\\
        \langle\calM\left(\bq;\balpha\right)\bhw,\cdot\rangle&0
    \end{bmatrix}
    \begin{bmatrix}
        \bhq^{\dagger{}} \\ g^*
    \end{bmatrix}
    &=
    \begin{bmatrix}
        0 \\ 1
    \end{bmatrix},
\end{align}
\end{subequations}
where convergence is most robust when the variables $\bhv$ and $\bhw$ are good approximations of the critical direct and adjoint eigenmode pair. Manipulation of \cref{eq:minaug} then yields the definitions for $g$, $\bhq$, and $\bhq^\dagger$ used in the actual implementation,
\begin{subequations}\label{eq:bifaug}
\begin{align}\label{eq:bifauga}
    \frac{1}{g}&=\left\langle\calM{}\left(\bq;\balpha\right)\bhv,\left[\mathrm{i}\omega\calM\left(\bq;\balpha\right)+\partial_{\bq}\calR|_{\bq,\balpha}\left(\cdot\right)\right]^{-1}\calM{}\left(\bq;\balpha\right)\bhw\right\rangle,\\
    \left[\mathrm{i}\omega\calM\left(\bq;\balpha\right)+\partial_{\bq}\calR|_{\bq,\balpha}\left(\cdot\right)\right]\bhq&=g\calM{}\left(\bq;\balpha\right)\bhw{},\\
    \left[\mathrm{i}\omega\calM\left(\bq;\balpha\right)+\partial_{\bq}\calR|_{\bq,\balpha}\left(\cdot\right)\right]^{\dagger}\bhq^{\dagger}&=g^*\calM{}\left(\bq;\balpha\right)\bhv{}.
\end{align}
\end{subequations}
Note that only the action of the inverse of the linear operator is required in \cref{eq:bifauga}; it is never actually constructed. This procedure then allows bifurcation points to be identified from a Newton-like sequence of corrections based on the following minimally augmented system in real arithmetic,
\begin{subequations}\label{eq:bifpoint}
\begin{align}\label{eq:bifpointa}
\begin{bmatrix}
    \partial_{\bq}\calR|_{\bq_k,\balpha_k}\left(\cdot\right)&\partial_{\alpha}\calR|_{\bq_k,\balpha_k}&0\\
    \left\langle\Re\left\{\partial_{\bq}g|_{\bq_k,\balpha_k}\right\},\cdot\right\rangle&\Re\left\{\partial_{\alpha}g|_{\bq_k,\balpha_k}\right\}&\Re\left\{\partial_{\omega}g|_{\bq_k,\balpha_k}\right\}\\
    -\left\langle\Im\left\{\partial_{\bq}g|_{\bq_k,\balpha_k}\right\},\cdot\right\rangle&\Im\left\{\partial_{\alpha}g|_{\bq_k,\balpha_k}\right\}&\Im\left\{\partial_{\omega}g|_{\bq_k,\balpha_k}\right\}
\end{bmatrix}\begin{bmatrix}
    \bdq_k\\
    \delta\alpha_k\\
    \delta\omega_k
\end{bmatrix}&=\begin{bmatrix}
    \calR\left(\bq_k;\balpha_k\right)\\
    \Re\left\{g\right\}\\
    \Im\left\{g\right\}
\end{bmatrix},\\
\begin{bmatrix}
    \bq_{k+1}\\
    \alpha_{k+1}\\
    \omega_{k+1}
\end{bmatrix}=\begin{bmatrix}
    \bq_k\\
    \alpha_k\\
    \omega_k
\end{bmatrix}-\beff{}&\left(\begin{bmatrix}
    \bdq_k\\
    \delta\alpha_k\\
    \delta\omega_k
\end{bmatrix}\right),\qquad{}k=0,1,\hdots.
\end{align}
\end{subequations}
To determine the entries in the bottom rows of the augmented Jacobian above, manipulation of \cref{eq:criticalityvariable} and \cref{eq:minaug} can be used to show that:
\begin{subequations}\label{eq:augjacobian}
\begin{align}
    \partial_{\bq}g|_{\bq,\balpha}&=\left\langle\bhq{}^{\dagger{}},\left[\mathrm{i}\omega\partial_{\bq}\calM|_{\bq,\balpha}\left(\cdot\right)\bhq+\partial_{\bq\bq}\calR|_{\bq,\balpha}\left(\cdot,\bhq\right)\right]\right\rangle,\\
    \partial_{\alpha}g|_{\bq,\balpha}&=\left\langle\bhq{}^{\dagger{}},\mathrm{i}\omega\partial_{\alpha}\calM|_{\bq,\balpha}\bhq+\partial_{\alpha\bq}\calR|_{\bq,\balpha}\left(\bhq\right)\right\rangle,\\
    \partial_{\omega}g_{\bq,\balpha}&=\left\langle\bhq{}^{\dagger{}},\mathrm{i}\calM\left(\bq;\balpha\right)\bhq\right\rangle.
\end{align}
\end{subequations}
\Cref{eq:augjacobian} allow the augmented Jacobian operator in \cref{eq:bifpointa} to be written explicitly. Thus, in \texttt{ff-bifbox}, \cref{eq:bifpoint} is solved using a Newton-like method leveraging \texttt{SNES} with an exact block Schur decomposition via \texttt{PCFIELDSPLIT} similar to the approach described in \cref{ssec:fixedpoints} for \cref{eq:MPcorrector}. Though the (unaugmented) Jacobian may become highly ill-conditioned as the solver converges to a bifurcation point, the basic Schur factorization is typically sufficient for Newton-like iteration techniques, and modified factorization schemes with greater robustness are also available \cite[Ch. 3.6]{govaerts_2000}. Once converged, the critical modes are normalized to the conditions $\sqrt{\left\langle\bhq,\calM_{\bq,\balpha}\bhq\right\rangle}=1$ and $\left\langle\bhq^\dagger,\calM_{\bq,\balpha}\bhq\right\rangle=1$.

As noted previously, the above approach is sufficient to describe any codimension-1 local bifurcation of a dynamical system of the form \cref{eq:nonlinearPDE}. However, for saddle--node and pitchfork bifurcations, the approach may be simplified since $\omega=0$. In these cases, \texttt{ff-bifbox} eliminates the final row of \cref{eq:bifpoint} (and the final column of the Jacobian operator), reducing the complexity of the block Schur factorization required from \texttt{PCFIELDSPLIT}. Additionally, the minimally augmented system \cref{eq:bifpoint} can be extended by defining additional criticality variables associated with each critical eigenmode and appending further augmentations to the augmented state vector, allowing handling of local bifurcations with higher codimension. Finally, just as fixed-point solutions satisfying \cref{eq:equilibrium} may be continued along a parameter using the predictor--corrector scheme based on \cref{eq:MPtangent} and \cref{eq:MPcorrector}, \texttt{ff-bifbox} supports tracing of bifurcations through the parameter space using the predictor--corrector scheme from \cref{sssec:tracing} with the state vector, Jacobian, and residual terms in \cref{eq:MPcorrector} replaced by their counterparts from the minimally augmented system \cref{eq:bifpoint}.

Once a bifurcation point is identified, the nonlinear dynamics in its immediate vicinity can be rigorously reduced to an equivalent normal form of the bifurcation with relatively modest additional effort. This procedure, which represents a projection of the local nonlinear dynamics onto the center manifold, allows classification of nondegenerate bifurcations as supercritical or subcritical and the construction of so-called ``weakly nonlinear'' solutions. Following~\cite{govaerts_2000}, these center manifold amplitude equations are written in their Poincar\'e normal forms as,
\begin{align}\label{eq:foldform}
    \frac{dA}{dt}+\partial_{\balpha}\lambda\bcdot\Delta{}\balpha+\beta{}A^2&=0,\\
    \frac{dA}{dt}+A\left(\partial_{\balpha}\lambda\bcdot\Delta{}\balpha+\beta{}\left|A\right|^2\right)&=0,
\end{align}
respectively, for a quadratic form (e.g. at a saddle--node bifurcation) and a cubic form (e.g. at a Hopf or pitchfork bifurcation). Here, $A$ is the (possibly complex) amplitude, $\partial_{\balpha}\lambda$ is a vector representing the gradient of the critical eigenvalue in the parameter space, $\Delta{}\balpha$ is the vector displacement in parameter space between a given point $\balpha$ and the critical point, and $\beta$ is the nonlinear coefficient. Making use of the above normalization conditions, these coefficients are calculated as,
\begin{subequations}\label{eq:formcoeffs}
\begin{align}
    \partial_{\balpha}\lambda&=\left\langle\bhq{}^\dagger{},\partial_{\balpha}\calR|_{\bq,\balpha}\right\rangle,\\
    \beta&=\frac{1}{2}\left\langle\bhq{}^\dagger{},\partial_{\bq\bq}\calR|_{\bq,\balpha}\left(\bhq{},\bhq\right)\right\rangle,
\end{align}
\end{subequations}
for the quadratic form, and,
\begin{subequations}\label{eq:formcoeffs2}
\begin{align}
    \partial_{\balpha}\lambda&=\left\langle\bhq{}^\dagger{},\mathrm{i}\omega\,\partial_{\bq}\calM|_{\bq,\balpha}\left(\bhq_{\alpha}\right)\bhq+\partial_{\bq\bq}\calR|_{\bq,\balpha}\left(\bhq_{\alpha},\bhq{}\right)\right\rangle+\left\langle\bhq{}^\dagger{},\mathrm{i}\omega\,\partial_{\balpha}\calM|_{\bq,\balpha}\,\bhq+\partial_{\balpha\bq}\calR|_{\bq,\balpha}\left(\bhq\right)\right\rangle,\\
    \beta&=\frac{1}{2}\left\langle\bhq{}^\dagger{},\mathrm{i}\omega\,\partial_{\bq\bq}\calM|_{\bq,\balpha}\left(\bhq{},\bhq^*\right)\bhq+\partial_{\bq\bq\bq}\calR|_{\bq,\balpha}\left(\bhq{},\bhq^*,\bhq\right)\right\rangle\nonumber\\
    &\hspace{1.5cm}+\frac{1}{2}\left\langle\bhq{}^\dagger{},\mathrm{i}\omega\,\partial_{\bq}\calM|_{\bq,\balpha}\left(\bhq_{\hat{q}^2}\right)\bhq{}^*+\partial_{\bq\bq}\calR|_{\bq,\balpha}\left(\bhq_{\hat{q}^2},\bhq{}^*\right)\right\rangle+\left\langle\bhq{}^\dagger{},\partial_{\bq\bq}\calR|_{\bq,\balpha}\left(\bhq{},\bhq_{|\hat{q}|^2}\right)\right\rangle,
\end{align}
for the cubic form, where the second-order corrections (cf. \cite{sipp_lebedev_2007}) to the bifurcating mode are,
\begin{align}
    \bhq_{\balpha}&=\left[\partial_{\bq}\calR|_{\bq,\balpha}\right]^{-1}\partial_{\balpha}\calR|_{\bq,\balpha},\\
    \bhq_{\hat{q}^2}&=\left[2\mathrm{i}\omega\calM|_{\bq,\balpha}+\partial_{\bq}\calR|_{\bq,\balpha}\right]^{-1}\left[\mathrm{i}\omega\,\partial_{\bq}\calM|_{\bq,\balpha}\left(\bhq\right)\bhq+\partial_{\bq\bq}\calR|_{\bq,\balpha}\left(\bhq,\bhq\right)\right],\\
    \bhq_{|\hat{q}|^2}&=\left[\partial_{\bq}\calR|_{\bq,\balpha}\right]^{-1}\partial_{\bq\bq}\calR|_{\bq,\balpha}\left(\bhq,\bhq^*\right).
\end{align}
\end{subequations}
Note that, in the cubic form, the real part of $\beta$ is the first Lyapunov coefficient, which indicates the super- or subcritical nature of the bifurcation when $\Re\{\beta\}>0$ or $\Re\{\beta\}<0$, respectively. The implementation in \texttt{ff-bifbox} automatically calculates the normal form associated with each bifurcation, allowing straightforward detection of codimension-2 cusp, Bautin, fold--Hopf, and Bogdanov--Takens bifurcations based on the behavior of the normal form coefficients during continuation. Though omitted here for brevity, normal form reductions are also implemented in \texttt{ff-bifbox} for fold--Hopf and Hopf--Hopf bifurcations. 
%%%%%%%%%%%%%%%%%%%%%%%%%%%
\subsection{Periodic solutions}

A Hopf bifurcation signals the intersection of branches of fixed-point and periodic solutions. Robust tracing along such (possibly unstable) periodic solution branches requires a representation of \cref{eq:nonlinearPDE} and its solutions in a manner that is constrained to be periodic. While multiple expansion bases and numerical techniques can be used to accomplish this (see~\cite{kuznetsov_2023} for an overview) and may be incorporated in future developments, \texttt{ff-bifbox} currently employs the harmonic balance method, leveraging a Fourier--Galerkin expansion of the system and state vector. This separation of variables allows the state vector to be expressed by the discrete Fourier series,
\begin{align}\label{eq:HBexpansion}
    \bq\left(\bx,t\right)&=\bhq_0\left(\bx\right)+\sum_{n=1}^{\infty}\left[\bhq_n\left(\bx\right)\exp\left(\mathrm{i}n\omega{}t\right)+\bhq_n^*\left(\bx\right)\exp\left(-\mathrm{i}n\omega{}t\right)\right],
\end{align}
where $\omega$ is the (unknown) fundamental angular frequency of the periodic solution. Note that \cref{eq:HBexpansion} is exact for a $T=2\pi/\omega$-periodic state.

Substituting the expansion \cref{eq:HBexpansion} into \cref{eq:nonlinearPDE}, a similar Fourier--Galerkin expansion for \cref{eq:nonlinearPDE} is found,
\begin{align}\label{eq:HBPDE}
    \calF{}\left(\bq,\bqd;\balpha\right)=\widehat{\calF}_0+\sum_{n=1}^{\infty}\left[\widehat{\calF}_n\exp\left(\mathrm{i}n\omega{}t\right)+\widehat{\calF}_n^*\exp\left(-\mathrm{i}n\omega{}t\right)\right]=0,
\end{align}
where $\widehat{\calF}_n=\widehat{(\calM{}\bqd)}_n+\widehat{\calR}_n$ denotes the $n$th harmonic component of the Fourier expansion of \cref{eq:nonlinearPDE}. Importantly, \cref{eq:HBPDE} is exact for a period-one solution, just as \cref{eq:equilibrium} is exact for a steady solution. Further, if $\calF{}$ is restricted to be analytic over the domain of the periodic orbit $\bq$, then the Fourier expansion of $\calF{}$ can be exactly constructed from the Taylor expansion of $\calF{}$ about its mean state $\bhq_0$,
\begin{align}\label{eq:HBF}
    &\calF{}\left(\bq,\bqd;\balpha\right)=\frac{1}{0!}\calR{}\left(\bhq_0;\balpha\right)\nonumber\\
    &\quad+\frac{1}{1!}\sum_{j=1}^{\infty}\left\{\left[\mathrm{i}j\omega\calM\left(\bhq_0;\balpha\right)\bhq_j+\partial_{\bq}\calR{}|_{\bhq_0,\balpha}\left(\bhq_j\right)\right]e^{\mathrm{i}j\omega{}t}+\left[-\mathrm{i}j\omega\calM\left(\bhq_0;\balpha\right)\bhq_j^*+\partial_{\bq}\calR{}^*|_{\bhq_0,\balpha}\left(\bhq_j^*\right)\right]e^{-\mathrm{i}j\omega{}t}\right\}\nonumber\\
    &\quad+\frac{1}{2!}\sum_{j=1}^{\infty}\sum_{k=1}^{\infty}\left\{\left[\mathrm{i}\left(j+k\right)\omega\partial_{\bq}\calM|_{\bhq_0,\balpha}\left(\bhq_k\right)\bhq_j+\partial_{\bq\bq}\calR{}_{\bhq_0,\balpha}\left(\bhq_j,\bhq_k\right)\right]e^{\mathrm{i}(j+k)\omega{}t}\right.\nonumber\\
    &\qquad\qquad\qquad+\left.\left[\mathrm{i}\left(j-k\right)\omega\partial_{\bq}\calM|_{\bhq_0,\balpha}\left(\bhq_k^*\right)\bhq_j+\partial_{\bq\bq}\calR{}_{\bhq_0,\balpha}\left(\bhq_j,\bhq_k^*\right)\right]e^{\mathrm{i}(j-k)\omega{}t}\right.\nonumber\\
    &\qquad\qquad\qquad+\left.\left[\mathrm{i}\left(-j+k\right)\omega\partial_{\bq}\calM|_{\bhq_0,\balpha}\left(\bhq_k\right)\bhq_j^*+\partial_{\bq\bq}\calR{}_{\bhq_0,\balpha}\left(\bhq_j^*,\bhq_k\right)\right]e^{\mathrm{i}(-j+k)\omega{}t}\right.\nonumber\\
    &\qquad\qquad\qquad+\left.\left[-\mathrm{i}\left(j+k\right)\omega\partial_{\bq}\calM|_{\bhq_0,\balpha}\left(\bhq_k^*\right)\bhq_j^*+\partial_{\bq\bq}\calR{}_{\bhq_0,\balpha}\left(\bhq_j^*,\bhq_k^*\right)\right]e^{-\mathrm{i}(j+k)\omega{}t}\right\}\nonumber\\
    &\quad+\frac{1}{3!}\sum_{j=1}^{\infty}\sum_{k=1}^{\infty}\sum_{l=1}^{\infty}\left\{\left[\mathrm{i}\left(j+k+l\right)\omega\partial_{\bq\bq}\calM|_{\bhq_0,\balpha}\left(\bhq_k,\bhq_l\right)\bhq_j+\partial_{\bq\bq\bq}\calR{}_{\bhq_0,\balpha}\left(\bhq_j,\bhq_k,\bhq_l\right)\right]e^{\mathrm{i}(j+k+l)\omega{}t}\right.\nonumber\\
    &\qquad\qquad\qquad{}\qquad{}+\left.\left[\mathrm{i}\left(j+k-l\right)\omega\partial_{\bq\bq}\calM|_{\bhq_0,\balpha}\left(\bhq_k,\bhq_l^*\right)\bhq_j+\partial_{\bq\bq\bq}\calR{}_{\bhq_0,\balpha}\left(\bhq_j,\bhq_k,\bhq_l^*\right)\right]e^{\mathrm{i}(j+k-l)\omega{}t}\right.\nonumber\\
    &\qquad\qquad{}\qquad{}\qquad{}+\left.\left[\mathrm{i}\left(j-k+k\right)\omega\partial_{\bq\bq}\calM|_{\bhq_0,\balpha}\left(\bhq_k^*,\bhq_l\right)\bhq_j+\partial_{\bq\bq\bq}\calR{}_{\bhq_0,\balpha}\left(\bhq_j,\bhq_k^*,\bhq_l\right)\right]e^{\mathrm{i}(j-k+l)\omega{}t}\right.\nonumber\\
    &\qquad\qquad\qquad{}\qquad{}+\left.\left[\mathrm{i}\left(-j+k+l\right)\omega\partial_{\bq\bq}\calM|_{\bhq_0,\balpha}\left(\bhq_k,\bhq_l\right)\bhq_j^*+\partial_{\bq\bq\bq}\calR{}_{\bhq_0,\balpha}\left(\bhq_j^*,\bhq_k,\bhq_l\right)\right]e^{\mathrm{i}(-j+k+l)\omega{}t}\right.\nonumber\\
    &\qquad\qquad\qquad{}\qquad{}+\left.\left[\mathrm{i}\left(-j+k-l\right)\omega\partial_{\bq\bq}\calM|_{\bhq_0,\balpha}\left(\bhq_k,\bhq_l^*\right)\bhq_j^*+\partial_{\bq\bq\bq}\calR{}_{\bhq_0,\balpha}\left(\bhq_j^*,\bhq_k,\bhq_l^*\right)\right]e^{\mathrm{i}(-j+k-l)\omega{}t}\right.\nonumber\\
    &\qquad\qquad\qquad{}\qquad{}+\left.\left[\mathrm{i}\left(-j-k+l\right)\omega\partial_{\bq\bq}\calM|_{\bhq_0,\balpha}\left(\bhq_k^*,\bhq_l\right)\bhq_j^*+\partial_{\bq\bq\bq}\calR{}_{\bhq_0,\balpha}\left(\bhq_j^*,\bhq_k^*,\bhq_l\right)\right]e^{\mathrm{i}(-j-k+l)\omega{}t}\right.\nonumber\\
    &\qquad\qquad\qquad{}\qquad{}+\left.\left[\mathrm{i}\left(j-k-l\right)\omega\partial_{\bq\bq}\calM|_{\bhq_0,\balpha}\left(\bhq_k^*,\bhq_l^*\right)\bhq_j+\partial_{\bq\bq\bq}\calR{}_{\bhq_0,\balpha}\left(\bhq_j,\bhq_k^*,\bhq_l^*\right)\right]e^{\mathrm{i}(j-k-l)\omega{}t}\right.\nonumber\\
    &\qquad\qquad\qquad{}\qquad{}+\left.\left[-\mathrm{i}\left(j+k+l\right)\omega\partial_{\bq\bq}\calM|_{\bhq_0,\balpha}\left(\bhq_k^*,\bhq_l^*\right)\bhq_j^*+\partial_{\bq\bq\bq}\calR{}_{\bhq_0,\balpha}\left(\bhq_j^*,\bhq_k^*,\bhq_l^*\right)\right]e^{\mathrm{i}(-j-k-l)\omega{}t}\right\}\nonumber\\
    &\quad{}+\hdots{}.
\end{align}
In the implementation, \texttt{ff-bifbox} truncates the continued Taylor series expansion in \cref{eq:HBF} at third order, limiting the exact periodic representation of $\calF{}$ to systems with polynomial nonlinearities of degree three or lower. However, so long as $\calF{}$ is analytic, auxiliary variables may be introduced to reduce the nonlinear system to cubic (or even quadratic) order, so this restriction does not necessitate approximation at the continuous level~\cite{cochelin_vergez_2009}.

Exploiting the orthogonality of the Fourier basis, a sequence of coupled nonlinear equations governing the evolution of each harmonic constituent of the periodic solution is obtained from \cref{eq:HBF} as $\widehat{\calF}_n=0$ for $n=0,1,2,\hdots$. More specifically, for the mean component ($n=0$) and unsteady tones ($n\geq1$), this procedure yields, respectively,
\begin{subequations}\label{eq:HBperiodic}
\begin{multline}\label{eq:HBperiodica}
    \calR{}\left(\bhq_0;\balpha\right)+\sum_{j=1}^{\infty}\left[\partial_{\bq\bq}\calR{}_{\bhq_0,\balpha}\left(\bhq_j,\bhq_j^*\right)\right]\\
    +\frac{1}{2}\sum_{j=1}^{\infty}\sum_{k=1}^{\infty}\left[\partial_{\bq\bq\bq}\calR{}_{\bhq_0,\balpha}\left(\bhq_j,\bhq_k,\bhq_{j+k}^*\right)+\partial_{\bq\bq\bq}\calR{}_{\bhq_0,\balpha}\left(\bhq_j^*,\bhq_k^*,\bhq_{j+k}\right)\right]=0,
\end{multline}
and
\begin{multline}\label{eq:HBperiodicb}
    \mathrm{i}n\omega\calM\left(\bhq_0;\balpha\right)\bhq_n+\partial_{\bq}\calR{}|_{\bhq_0,\balpha}\left(\bhq_n\right)\\
    +\frac{1}{2}\sum_{j=1}^{n-1}\left[\mathrm{i}n\omega\partial_{\bq}\calM|_{\bhq_0,\balpha}\left(\bhq_{j}\right)\bhq_{n-j}+\partial_{\bq\bq}\calR{}_{\bhq_0,\balpha}\left(\bhq_{j},\bhq_{n-j}\right)\right]
    +\sum_{j=n+1}^{\infty}\left[\partial_{\bq\bq}\calR{}_{\bhq_0,\balpha}\left(\bhq_{j},\bhq_{j-n}^*\right)\right]\\
    +\frac{1}{6}\sum_{j=1}^{\infty}\sum_{k=1}^{n-j-1}\left[\mathrm{i}n\omega\partial_{\bq\bq}\calM|_{\bhq_0,\balpha}\left(\bhq_j,\bhq_k\right)\bhq_{n-j-k}+\partial_{\bq\bq\bq}\calR{}_{\bhq_0,\balpha}\left(\bhq_j,\bhq_k,\bhq_{n-j-k}\right)\right]\\    
    +\frac{1}{2}\sum_{j=1}^{\infty}\sum_{k=j-n+1}^{\infty}\left[\mathrm{i}n\omega\partial_{\bq\bq}\calM|_{\bhq_0,\balpha}\left(\bhq_j,\bhq_k\right)\bhq_{j+k-n}^*+\partial_{\bq\bq\bq}\calR{}_{\bhq_0,\balpha}\left(\bhq_j,\bhq_k,\bhq_{j+k-n}^*\right)\right]\\
    +\frac{1}{2}\sum_{j=1}^{\infty}\sum_{k=1}^{j-n-1}\left[\mathrm{i}n\omega\partial_{\bq\bq}\calM|_{\bhq_0,\balpha}\left(\bhq_j,\bhq_k^*\right)\bhq_{j-k-n}^*+\partial_{\bq\bq\bq}\calR{}_{\bhq_0,\balpha}\left(\bhq_j,\bhq_k^*,\bhq_{j-k-n}^*\right)\right]=0.
\end{multline}
\Cref{eq:HBperiodica} and \cref{eq:HBperiodicb} describe a family of periodic solutions with indeterminate phase. To uniquely select a particular solution with fixed phase, \texttt{ff-bifbox} uses the integral phase constraint~\cite{kuznetsov_2023}, expressed in the frequency domain as,
\begin{align}\label{eq:HBperiodicc}
    \sum_{n=0}^{\infty}\left\langle\bhq{}_n,\partial_{\omega}\widehat{\calF}_{n}\right\rangle&=0.
\end{align}
\end{subequations}

The exact expansion \cref{eq:HBexpansion} and corresponding equation sequence \cref{eq:HBperiodic} can then be truncated at some finite $n=N$ to approximate the periodic state and bound the dimension of the system. Importantly, the exponential convergence property of the spectral method allows temporal convergence to be achieved with relatively modest values of $N$ compared to other approaches involving, for example, finite difference approximations of the temporal derivatives~\cite{boyd_2013, detroux_etal_2015, kuznetsov_2023}. Once truncated, the solution of \cref{eq:HBperiodic} follows a Newton-like method leveraging PETSc \texttt{SNES} in real arithmetic,
\begin{subequations}\label{eq:HBnewton}
\begin{align}
    \begin{bmatrix}
        \partial_{\bhq_0}\widehat{\calF}_{0}&\partial_{\Re\{\bhq_1\}}\widehat{\calF}_{0}&\partial_{\Im\{\bhq_1\}}\widehat{\calF}_{0}&\hdots&0\\
        \Re\{\partial_{\bhq_0}\widehat{\calF}_{1}\}&\Re\{\partial_{\Re\{\bhq_1\}}\widehat{\calF}_{1}\}&\Re\{\partial_{\Im\{\bhq_1\}}\widehat{\calF}_{1}\}&\hdots&\Re\{\partial_\omega\widehat{\calF}_{1}\}\\
        \Im\{\partial_{\bhq_0}\widehat{\calF}_{1}\}&\Im\{\partial_{\Re\{\bhq_1\}}\widehat{\calF}_{1}\}&\Im\{\partial_{\Im\{\bhq_1\}}\widehat{\calF}_{1}\}&\hdots&\Im\{\partial_\omega\widehat{\calF}_{1}\}\\
        \vdots&\vdots&\vdots&\ddots&\vdots&\\
        0&\left\langle\Re\{\partial_\omega\widehat{\calF}_{1}\},\cdot\right\rangle&\left\langle\Im\{\partial_\omega\widehat{\calF}_{1}\},\cdot\right\rangle&\hdots&0
    \end{bmatrix}_k\begin{bmatrix}
        \bdhq_{0}\\
        \Re\left\{\bdhq_{1}\right\}\\
        \Im\left\{\bdhq_{1}\right\}\\
        \vdots\\
        \delta\omega
    \end{bmatrix}_k&=\begin{bmatrix}
        \widehat{\calF}_{0}\\
        \Re\{\widehat{\calF}_{1}\}\\
        \Im\{\widehat{\calF}_{1}\}\\
        \vdots\\
        0
    \end{bmatrix}_k,\\
    \begin{bmatrix}
        \bhq_{0}\\
        \Re\left\{\bhq_{1}\right\}\\
        \Im\left\{\bhq_{1}\right\}\\
        \vdots\\
        \omega
    \end{bmatrix}_{k+1}=\begin{bmatrix}
        \bhq_{0}\\
        \Re\left\{\bhq_{1}\right\}\\
        \Im\left\{\bhq_{1}\right\}\\
        \vdots\\
        \omega
    \end{bmatrix}_k-\beff{}\left(\begin{bmatrix}
        \bdhq_{0}\\
        \Re\left\{\bdhq_{1}\right\}\\
        \Im\left\{\bdhq_{1}\right\}\\
        \vdots\\
        \delta\omega
    \end{bmatrix}_k\right)&,\qquad{}k=0,1,\hdots.
\end{align}
\end{subequations}
\Cref{eq:HBnewton} can also be augmented or modified analogously to \cref{eq:Newton} to support branch tracing and branch switching techniques, as described in \cref{sssec:tracing} and \cref{sssec:switching}, respectively.

Importantly, the harmonic balance Jacobian operator is $2N+1$ times larger than the steady Jacobian operator, indicating a significant increase in computational expense for periodic solutions compared to steady solutions. To help manage this increased difficulty, \texttt{ff-bifbox} interfaces with \texttt{PCFIELDSPLIT} from PETSc to support efficient iterative resolution approaches with block preconditioners. Some additional discussion related to these matters can be found in Rigas \textit{et al.}~\cite{rigas_etal_2021}, Sierra \textit{et al.}~\cite{sierra-ausin_etal_2022}, and prior work~\cite{douglas_etal_2021b}.

%%%%%%%%%%%%%%%%%%%%%%%%%%%
\subsection{Stability and bifurcations of periodic solutions}\label{ssec:stability_bifurcations_periodic}

The stability of $T$-periodic solutions to \cref{eq:nonlinearPDE} may be studied through the asymptotic time evolution of perturbations. As in \cref{ssec:stability}, any such perturbation may be written as a superposition of modes. Here, the perturbations are expressed generically in terms of $T$-periodic Floquet modes $\bbq$ of the form:
\begin{align}\label{eq:HBsuperposition}
    \baq\left(\bx,t\right)&=\sum\left[\bbq{}\left(\bx,t\right)\exp\left(\lambda{}t\right)+\bbq{}^*\left(\bx,t\right)\exp\left(\lambda{}^*t\right)\right].
\end{align}
Then, expanding the space-time modes $\bbq$ using \cref{eq:HBexpansion}, the arbitrary perturbation may be expressed as,
\begin{multline}\label{eq:HBsuperposition2}
    \baq\left(\bx,t\right)=\sum\left[\left\{\bhbq_0\left(\bx\right)+\sum_{n=1}^{\infty}\left[\bhbq_{n}\left(\bx\right)\exp\left(\mathrm{i}n\omega{}t\right)+\bhbq_{n}^*\left(\bx\right)\exp\left(-\mathrm{i}n\omega{}t\right)\right]\right\}\exp\left(\lambda{}t\right).\right.\\
    +\left.\left\{\bhbq_0^*\left(\bx\right)+\sum_{n=1}^{\infty}\left[\bhbq_{n}^*\left(\bx\right)\exp\left(-\mathrm{i}n\omega{}t\right)+\bhbq_{n}\left(\bx\right)\exp\left(\mathrm{i}n\omega{}t\right)\right]\right\}\exp\left(\lambda{}^*t\right)\right].
\end{multline}
Once again, the modes $\bbq$ are generally coupled and non-orthogonal, but they become decoupled when the perturbations are infinitesimally small, and, under linearity, they become orthogonal at asymptotically large times. This gives rise to the generalized eigenvalue problem governing the Floquet stability of the periodic orbit using Hill's method~\cite{hill_1886, detroux_etal_2015},
\begin{multline}\label{eq:HBfloquet}
    \left\langle\begin{bmatrix}
        \bhbq_0^\dagger\\\bhbq_1^\dagger\\\bhbq_1^{*\dagger}\\\vdots
    \end{bmatrix},\left[\lambda_F\begin{bmatrix}
        \partial_{\bhbq_0}\widehat{(\calM\bbq)}_{0}&\partial_{\bhbq_1}\widehat{(\calM\bbq)}_{0}&\partial_{\bhbq_1^*}\widehat{(\calM\bbq)}_{0}&\hdots\\
        \partial_{\bhbq_0}\widehat{(\calM\bbq)}_{1}&\partial_{\bhbq_1}\widehat{(\calM\bbq)}_{1}&\partial_{\bhbq_1^*}\widehat{(\calM\bbq)}_{1}&\hdots\\
        \partial_{\bhbq_0}\widehat{(\calM\bbq)}_{1}^*&\partial_{\bhbq_1}\widehat{(\calM\bbq)}_{1}^*&\partial_{\bhbq_1^*}\widehat{(\calM\bbq)}_{1}^*&\hdots\\
        \vdots&\vdots&\vdots&\ddots
    \end{bmatrix}+\begin{bmatrix}
        \partial_{\bhbq_0}\widehat{\calF}_{0}&\partial_{\bhbq_1}\widehat{\calF}_{0}&\partial_{\bhbq_1^*}\widehat{\calF}_{0}&\hdots\\
        \partial_{\bhbq_0}\widehat{\calF}_{1}&\partial_{\bhbq_1}\widehat{\calF}_{1}&\partial_{\bhbq_1^*}\widehat{\calF}_{1}&\hdots\\
        \partial_{\bhbq_0}\widehat{\calF}_{1}^*&\partial_{\bhbq_1}\widehat{\calF}_{1}^*&\partial_{\bhbq_1^*}\widehat{\calF}_{1}^*&\hdots\\
        \vdots&\vdots&\vdots&\ddots
    \end{bmatrix}\right]\begin{bmatrix}
        \bhbq_0\\\bhbq_1\\\bhbq_1^*\\\vdots
    \end{bmatrix}\right\rangle\\
    =\left\langle\bbq^\dagger,\left[\lambda_F\partial_{\bq}\left[\calM\left(\bq;\balpha\right)\bq\right]+\partial_{\bq}\calF|_{\bq,\balpha}\left(\cdot\right)\right]\bbq\right\rangle=0,
\end{multline}
where, $\bbq$ and $\bbq^{\dagger}$ represent the respective direct (right) and adjoint (left) Floquet eigenmode pair corresponding to a particular eigenvalue (i.e. Floquet exponent) $\lambda_F=\sigma_F+\mathrm{i}\omega_F$ with growth rate $\Re\left\{\lambda\right\}=\sigma_F$ and frequency $\Im\left\{\lambda\right\}=\omega_F$.

Note that, at the continuous level, Floquet exponents appear as an infinity of equivalent eigenvalues with frequencies shifted by an integer number of periods $\lambda_F=\sigma_F+\mathrm{i}\left(\omega_F+k\omega\right), k\in\mathbb{Z}$. For example, the periodic solution about which stability is computed is associated with an infinity of neutrally stable Floquet exponents $\lambda_F=\mathrm{i}k\omega$. However, the truncation of the invariant periodic state to a finite number of harmonics destroys this invariance, causing eigenvalues with different integer-period shifts to have a finite number of non-identical Floquet exponents. For more discussion about the interpretation of Floquet analysis using Hill's method, the reader is referred to Moulin~\cite{moulin_2020} and Sierra \textit{et al.}~\cite{sierra_etal_2021}.

Once an instability or near-instability of a periodic solution is identified using Floquet analysis, the exact state-space location and dynamical properties of the nearby local bifurcation of the periodic orbit may be obtained by enforcing a criticality condition. The procedure for computing and continuing bifurcations of periodic orbits is essentially identical to the corresponding process for bifurcations of fixed-points discussed in \cref{ssec:bifurcations} with the steady operators replaced by periodic operators. This further allows for normal form reduction and center manifold projection of the nonlinear state in the vicinity of the bifurcations.

%%%%%%%%%%%%%%%%%%%%%%%%%%%
\section{Validation and examples} \label{sec:validation}
%%%%%%%%%%%%%%%%%%%%%%%%%%%

The \texttt{ff-bifbox} package is extensively validated against a growing collection of prior investigations~\cite{douglas_etal_2021b,douglas_etal_2022,meliga_etal_2012,sipp_lebedev_2007,brokof_etal_2024,chevalier_etal_2024,douglas_2024,douglas_etal_2023,fani_etal_2018,garnaud_2012,marquet_larsson_2015,moulin_etal_2019,wang_etal_2024,pralits_etal_2010,sipp_marquet_2013,moreno-boza_etal_2018}. These examples focus on two- and three-dimensional fluid mechanics problems spanning incompressible, compressible, and reacting flow. Working examples/tutorials capable of reproducing the main results of each of these studies is currently available on the \texttt{ff-bifbox} repository at \href{github.com/cmdoug/ff-bifbox}{github.com/cmdoug/ff-bifbox}, and additional examples from literature are added regularly.

In this work, exemplary results are provided for three nonlinear PDE systems: the ``Brusselator'' reaction--diffusion system, an elastic system with geometric nonlinearity, and a compressible Navier--Stokes system. The \texttt{ff-bifbox} configuration files, FreeFEM mesh generation scripts, and additional documentation for each case are available on the repository at \href{github.com/cmdoug/ff-bifbox}{github.com/cmdoug/ff-bifbox}. While the library is designed to be scalable for large problems on distributed-memory systems, all results presented here have been computed on a laptop computer with 48 GB of RAM and 1$\times$14-core CPU (Apple M4 Pro) to ensure accessibility.

It is worth noting that the main computational cost of the schemes described in \cref{sec:theory} arises from the resolution of large, sparse, and distributed linear systems. While \texttt{PCFIELDSPLIT} provides a convenient interface to parameterize inner block preconditioners and solvers that is accessible in \texttt{ff-bifbox}, the design of efficient nested block methods lies beyond the scope of this communication. In the following results, attention is therefore restricted to exact block solvers, namely \texttt{PCLU} and \texttt{PCCHOLESKY}. With the advent of advanced techniques such as block low-rank multifrontal factorizations~\cite{amestoy_etal_2017}---readily available in PETSc---these block solvers remain well-suited for tackling large-scale systems. This is demonstrated for the solution and tracing of periodic orbits in \cref{ssec:cylinder}.

%%%%%%%%%%%%%%%%%%%%%%%%%%%
\subsection{Dynamics of the Brusselator system in three dimensions}\label{ssec:brusselator}

The Brusselator system is a nonlinear PDE modeling an autocatalytic chemical reaction--diffusion system originally formulated by Prigogine and Lefever~\cite{prigogine_lefever_1968}. It is a well-studied example of non-equilibrium chemical dynamics and pattern formation~\cite{holodniok_etal_1987,lust_roose_2000,pena_perez-garcia_2001,uecker_2021} that admits partial analytic treatment on elementary domains, making it a useful validation case. The system may be written as,
\begin{subequations}\label{eq:brusselator}
\begin{align}
    \dot{X}&=A-\left(B+1\right)X+X^2Y+\frac{D_X}{L^2}\nabla^2X\\
    \dot{Y}&=BX-X^2Y+\frac{D_Y}{L^2}\nabla^2Y,
\end{align}
\end{subequations}
where $\bq=\left[X,Y\right]^T$ is the state vector and $\balpha=\left[A,B,D_X,D_Y,L\right]^T$ are the positive parameters. Clearly, \cref{eq:brusselator} is of the form \cref{eq:nonlinearPDE}. Here, following Lust \& Roose~\cite{lust_roose_2000}, the parameters $A=2$, $B=5.45$, $D_X=0.008$, and $D_Y=0.004$ are fixed, while the length scale $L$ is varied. The system \cref{eq:brusselator} is posed in a cube with Dirichlet boundary conditions given by $\left(X,Y\right)=\left(A,B/A\right)$ on two opposite faces and homogeneous Neumann boundary conditions on the remaining faces. Importantly for validation, this 3-D configuration admits the entire subset of solutions spanned by the 1-D Brusselator investigated by Lust \& Roose~\cite{lust_roose_2000} (in addition to many other 2-D and 3-D solutions). For example, the system has as a ``base state'' the spatially-uniform equilibrium solution: $\left(X,Y\right)=\left(A,B/A\right)$. However, it will be verified that, as $L$ increases, this base state becomes destabilized, causing the system to spontaneously express temporal oscillations and spatial patterns.

Analytically, the system's leading bifurcations are determined from its characteristic polynomial derived via linearization of \cref{eq:brusselator} about the base state. Such bifurcations occur whenever an eigenvalue with zero real part satisfies,
\begin{align}\label{eq:charpoly}
    \left[\lambda+1-B+\frac{D_X}{L^2}\pi^2\left(k^2+l^2+m^2\right)\right]
    \left[\lambda+A^2+\frac{D_Y}{L^2}\pi^2\left(k^2+l^2+m^2\right)\right]+A^2B&=0,
\end{align}
where $k$ is any positive integer and $l,m$ are any non-negative integers. For example, a sequence of Hopf bifurcations of the base state occur at purely imaginary roots of \cref{eq:charpoly},
\begin{align}\label{eq:hopf}
    L&=\pi\sqrt{\frac{D_X+D_Y}{B-A^2-1}}\sqrt{k^2+l^2+m^2}\approx0.51302\sqrt{k^2+l^2+m^2},
\end{align}
where the associated $\Im\{\lambda\}=\omega=\sqrt{A^2B-\left[\frac{A^2D_X+\left(B-1\right)D_Y}{D_X+D_Y}\right]^2}\approx2.1395$. Similarly, steady pitchfork bifurcations of the base state would occur at zero roots of \cref{eq:charpoly},
\begin{align}\label{eq:pitchfork}
    \left[A^2+\frac{D_Y}{L^2}\pi^2\left(k^2+l^2+m^2\right)\right]\left[1+\frac{D_X}{L^2}\pi^2\left(k^2+l^2+m^2\right)\right]-B\frac{D_Y}{L^2}\pi^2\left(k^2+l^2+m^2\right)&=0.
\end{align}
However, it can be shown that no real $L$ can satisfy \cref{eq:pitchfork} under the chosen set of parameters.

A subset of the dynamics for the considered 3-D Brusselator system for $L\leq2$ are calculated using \texttt{ff-bifbox} as described in \cref{sec:theory} with discretization via FreeFEM~\cite{hecht_2012} based on piecewise quadratic tetrahedral finite elements in a cube with $11\times11\times11$ vertices. The square prismatic symmetry of the geometry and boundary conditions is exploited to simplify numerics, but all calculations are performed in 3-D. The step-by-step commands used to compute the base state, analyze stability, identify bifurcation points, trace periodic orbits, and characterize Floquet stability are provided in the public repository. These numerical results are then compared against the analytical Hopf points from \cref{eq:hopf}, as well as the 1-D numerical results of Lust \& Roose~\cite{lust_roose_2000}. Since the purpose of this example is to validate the solver, a comprehensive characterization of the full 3-D Brusselator state space is not sought. Instead, the presented results concentrate on tracing periodic solution branches for which validation data is available, i.e. those emanating from analytically-determined Hopf bifurcations of the base state or from the pitchfork bifurcations of the 1-D solution branches reported by Lust \& Roose~\cite{lust_roose_2000}. Floquet analysis is also used to determine assess the linear stability of the 1-D solutions, revealing additional pitchfork bifurcations to 2-D and 3-D periodic orbits as well as Neimark--Sacker (NS) bifurcations to 1-D, 2-D, and 3-D quasiperiodic states.

\begin{figure}
    \centering
    \input{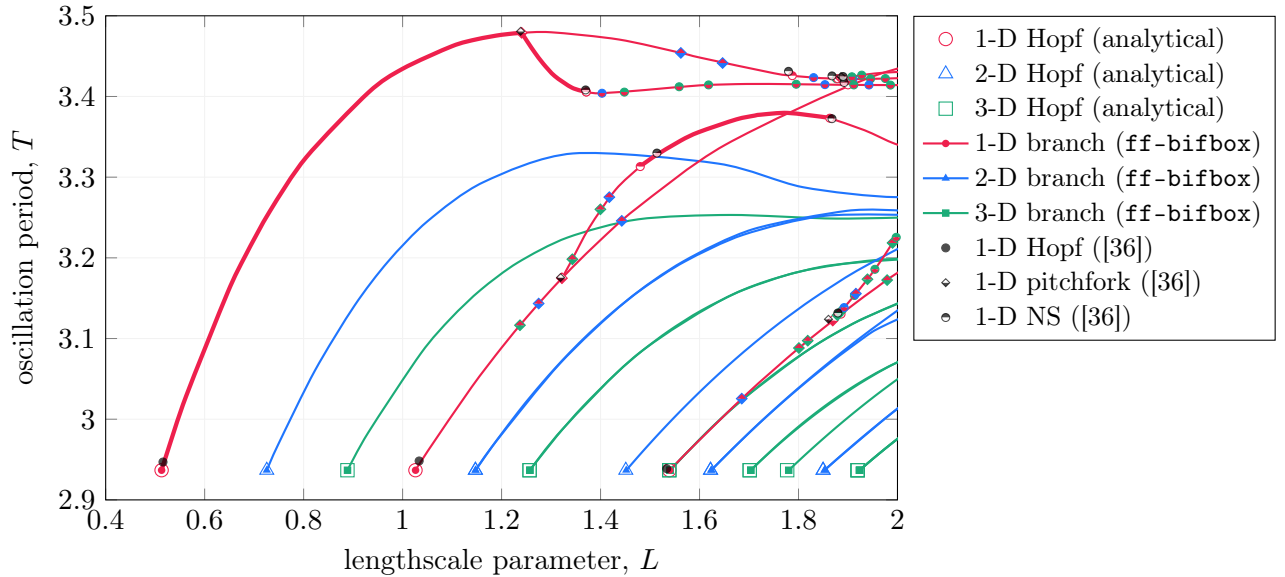}
    \caption{Bifurcation diagram showing the oscillation period $T=2\pi/\omega$ of selected time-periodic solutions to the 3-D Brusselator system with varying $L$ for $A=2$, $B=5.45$, $D_X=0.008$, and $D_Y=0.004$. Respectively, pitchfork and Neimark--Sacker (NS) bifurcations of periodic orbits are denoted using two-tone diamonds and circles; red, blue, and green lines and markers denote 1-D, 2-D, and 3-D branches; and thick and thin lines denote asymptotically stable and unstable solution branches. Hollow and solid markers along $T\approx2.94$ represent the exact and numerical Hopf points, respectively. Hopf, pitchfork, and NS bifurcations of the 1-D system as reported by Lust \& Roose~\cite{lust_roose_2000} are overlaid in black.}
    \label{fig:brusselator}
\end{figure}

A bifurcation diagram summarizing the computed dynamics of the 3-D Brusselator system is shown in \cref{fig:brusselator}. First, it may be observed that each of the Hopf bifurcations identified analytically is accurately captured using \texttt{ff-bifbox}, as indicated by the agreement between the hollow and solid markers along the bottom of \cref{fig:brusselator}, from which emanate various branches of 1-D, 2-D, and 3-D periodic solutions.

Next, the results from \texttt{ff-bifbox} may be compared against the 1-D solution branches originally characterized by Lust \& Roose~\cite{lust_roose_2000}. As in their work, seven distinct branches of 1-D solutions arise for $L\leq2$. All twelve of the 1-D pitchfork and NS bifurcations identified in that prior study are recovered using \texttt{ff-bifbox}, in addition to many additional pitchfork and NS bifurcations associated with 2-D or 3-D solutions. As evidenced by the increasing density of bifurcations (and, therefore, solution branches) with increasing $L$, the dynamics of the system become substantially more complex as the system becomes less dissipative.

%%%%%%%%%%%%%%%%%%%%%%%%%%%
\subsection{Finite-amplitude static deformation of a three-dimensional elastic structure}\label{ssec:structure}

\begin{figure}
    \centering
    \tdplotsetmaincoords{70}{130}
    \begin{tikzpicture}[tdplot_main_coords, scale=17,trig format=rad]
    % definitions
    \def\L{0.508};
    \def\R{2.540};
    \def\B{0.1};
    \def\t{0.00635};
    \def\labang{0.07};
    %% x=const faces
    \draw[thin,draw=black,fill=paperorange,nearly transparent] (-\L/2,{-(\R+\t/2)*sin(\B)},{\t/2*cos(\B)}) plot[variable=\x,domain=-\B:\B,smooth](-\L/2,{(\R+\t/2)*sin(\x)},{-\R*cos(\B)+(\R+\t/2)*cos(\x)}) -- (-\L/2,{(\R-\t/2)*sin(\B)},{-\t/2*cos(\B)}) plot[variable=\x,domain=\B:-\B,smooth](-\L/2,{(\R-\t/2)*sin(\x)},{-\R*cos(\B)+(\R-\t/2)*cos(\x)}) -- (-\L/2,{-(\R+\t/2)*sin(\B)},{\t/2*cos(\B)});
    \draw[draw=black,fill=paperorange,nearly transparent] (\L/2,{-(\R+\t/2)*sin(\B)},{\t/2*cos(\B)}) plot[variable=\x,domain=-\B:\B,smooth](\L/2,{(\R+\t/2)*sin(\x)},{-\R*cos(\B)+(\R+\t/2)*cos(\x)}) -- (\L/2,{(\R-\t/2)*sin(\B)},{-\t/2*cos(\B)}) plot[variable=\x,domain=\B:-\B,smooth](\L/2,{(\R-\t/2)*sin(\x)},{-\R*cos(\B)+(\R-\t/2)*cos(\x)}) -- (\L/2,{-(\R+\t/2)*sin(\B)},{\t/2*cos(\B)});
    %% angle=const faces
    \draw[draw=black,fill=paperorange,nearly transparent] (-\L/2,{-(\R+\t/2)*sin(\B)},{\t/2*cos(\B)}) -- (\L/2,{-(\R+\t/2)*sin(\B)},{\t/2*cos(\B)}) -- (\L/2,{-(\R-\t/2)*sin(\B)},{-\t/2*cos(\B)}) -- (-\L/2,{-(\R+\t/2)*sin(\B)},{-\t/2*cos(\B)}) -- cycle;
    \draw[draw=black,fill=paperorange,nearly transparent] (-\L/2,{(\R+\t/2)*sin(\B)},{\t/2*cos(\B)}) -- (\L/2,{(\R+\t/2)*sin(\B)},{\t/2*cos(\B)}) -- (\L/2,{(\R-\t/2)*sin(\B)},{-\t/2*cos(\B)}) -- (-\L/2,{(\R+\t/2)*sin(\B)},{-\t/2*cos(\B)}) -- cycle;
    %% radius=const faces
    \draw[draw=black,fill=paperorange,nearly transparent] (-\L/2,{-(\R+\t/2)*sin(\B)},{\t/2*cos(\B)}) plot[variable=\x,domain=-\B:\B,smooth](-\L/2,{(\R+\t/2)*sin(\x)},{-\R*cos(\B)+(\R+\t/2)*cos(\x)}) -- (\L/2,{(\R+\t/2)*sin(\B)},{\t/2*cos(\B)}) plot[variable=\x,domain=\B:-\B,smooth](\L/2,{(\R+\t/2)*sin(\x)},{-\R*cos(\B)+(\R+\t/2)*cos(\x)}) -- (-\L/2,{-(\R+\t/2)*sin(\B)},{\t/2*cos(\B)});
    \draw[draw=black,fill=paperorange,nearly transparent] (-\L/2,{-(\R-\t/2)*sin(\B)},{-\t/2*cos(\B)}) plot[variable=\x,domain=-\B:\B,smooth](-\L/2,{(\R-\t/2)*sin(\x)},{-\R*cos(\B)+(\R-\t/2)*cos(\x)}) -- (\L/2,{(\R-\t/2)*sin(\B)},{-\t/2*cos(\B)}) plot[variable=\x,domain=\B:-\B,smooth](\L/2,{(\R-\t/2)*sin(\x)},{-\R*cos(\B)+(\R-\t/2)*cos(\x)}) -- (-\L/2,{-(\R-\t/2)*sin(\B)},{-\t/2*cos(\B)});
    %%coordinate axes
    \draw[very thin, dotted] (0,0,0) -- (0.27,0,0);
    \draw[very thin, dotted] (0,0,0) -- (0,0.27,0);
    \draw[very thin, dotted] (0,0,0) -- (0,0,0.135);
    \draw[thick,-stealth] (0.27,0,0) --++ (0.075,0,0) node[below]{$x$};
    \draw[thick,-stealth] (0,0.27,0) --++ (0,0.075,0) node[right]{$y$};
    \draw[thick,-stealth] (0,0,0.135) --++ (0,0,0.075) node[left]{$z$};
    % point load
    \draw[ultra thick, -stealth] (0,0,0.1) node[left]{$P$} -- (0,0,{\t/2+\R*(1-cos(\B)});
    % lines of symmetry
    \draw[thin, dash dot] plot[variable=\x,domain=-\B:\B,smooth](0,{(\R+\t/2)*sin(\x)},{-\R*cos(\B)+(\R+\t/2)*cos(\x)}) --++ (0,0,-\t);
    \draw[thin, dash dot,nearly transparent] plot[variable=\x,domain=\B:-\B,smooth](0,{(\R-\t/2)*sin(\x)},{-\R*cos(\B)+(\R-\t/2)*cos(\x)}) --++ (0,0,\t);
    \draw[thin, dash dot] (\L/2,0,{-\t/2+\R*(1-cos(\B))}) --++ (0,0,\t) --++ (-\L,0,0);
    \draw[thin, dash dot,nearly transparent] (-\L/2,0,{\t/2+\R*(1-cos(\B))}) --++ (0,0,-\t) --++ (\L,0,0);
    %pinned BCs
    \draw[thick] (-\L/2,{\R*sin(\B)},0) --++ (\L,0,0);
    \draw[thick] (-\L/2,{-\R*sin(\B)},0) --++ (\L,0,0);
    \draw[semitransparent,fill=gray] (\L/2,{\R*sin(\B)},0) --++(0,0.01,0.01)--++(0,0,-0.02)--cycle;
    \fill[pattern color=black, pattern=north west lines] (\L/2,{\R*sin(\B)+0.01},0.012) --++ (0,0.01,0) --++ (0,0,-0.024) --++ (0,-0.01,0) -- cycle;
    \draw[semitransparent,fill=gray] (-\L/2,{\R*sin(\B)},0) --++(0,0.01,0.01)--++(0,0,-0.02)--cycle;
    \fill[pattern color=black, pattern=north west lines] (-\L/2,{\R*sin(\B)+0.01},0.012) --++ (0,0.01,0) --++ (0,0,-0.024) --++ (0,-0.01,0) -- cycle;
    \draw[semitransparent,fill=gray] (\L/2,{-\R*sin(\B)},0) --++(0,-0.01,0.01)--++(0,0,-0.02)--cycle;
    \fill[pattern color=black, pattern=north west lines] (\L/2,{-\R*sin(\B)-0.01},0.012) --++ (0,-0.01,0) --++ (0,0,-0.024) --++ (0,0.01,0) -- cycle;
    \draw[semitransparent,fill=gray] (-\L/2,{-\R*sin(\B)},0) --++(0,-0.01,0.01)--++(0,0,-0.02)--cycle;
    \fill[pattern color=black, pattern=north west lines] (-\L/2,{-\R*sin(\B)-0.01},0.012) --++ (0,-0.01,0) --++ (0,0,-0.024) --++ (0,0.01,0) -- cycle;
    % L-dimension
    \draw [thin,stealth-stealth] (-\L/2,{\R*sin(\labang)},{-\R*cos(\B)+(\R+\t/2)*cos(\labang)}) node[above]{$L$} --++ (\L,0,0);
    % R-dimension
    \draw [thin,-stealth] (\L/2,{-(\R+0.08)*sin(\labang)},{(0.08-\R)*cos(\B)+(\R+\t/2)*cos(-\labang)}) node[above]{$R$} -- (\L/2,{-\R*sin(\labang)},{-\R*cos(\B)+(\R+\t/2)*cos(-\labang)});
    % t-dimension
    \draw [thin,-stealth](\L/2,{-(\R+0.03)*sin(\labang/2)},{(0.03-\R)*cos(\B)+(\R+\t/2)*cos(-\labang/2)}) -- (\L/2,{-\R*sin(\labang/2)},{-\R*cos(\B)+(\R+\t/2)*cos(-\labang/2)});
    \draw [thin,-stealth](\L/2,{-(\R-0.03)*sin(\labang/2)},{(-0.03-\R)*cos(\B)+(\R-\t/2)*cos(-\labang/2)}) node[left]{$t$} -- (\L/2,{-\R*sin(\labang/2)},{-\R*cos(\B)+(\R-\t/2)*cos(-\labang/2)});
    % \beta-dimension
    \draw[thin, stealth-stealth] plot[variable=\x,domain=-\B:0,smooth](\L/4,{(\R+\t/2)*sin(\x)},{-\R*cos(\B)+(\R+\t/2)*cos(\x)});
    \node at (\L/4,{(\R+\t/2)*sin(-\labang)},{-\R*cos(\B)+(\R+\t/2)*cos(\labang)}) [above right]{$\beta$};
    % Other constants
    \node at (0,0.35,0.25) {$\begin{aligned}
    R&=2540\,\text{mm}\\
    L&=508\,\text{mm}\\
    t&=6.35\,\text{mm}\\
    \end{aligned}\quad\begin{aligned}
    \beta&=0.1\,\text{rad}\\
    E&=3.10275\,\text{GPa}\\
    \nu&=0.3
    \end{aligned}$};
    \end{tikzpicture}
    \caption{Schematic of the benchmark problem for the hinged cylindrical section under point loading parameterized by $P$ with lines of symmetry indicated. Note that the pinned boundary conditions are enforced on a line only along the neutral axis of the body~\cite{leahu-aluas_abed-meriam_2011}.}
    \label{fig:bodygeometry}
\end{figure}
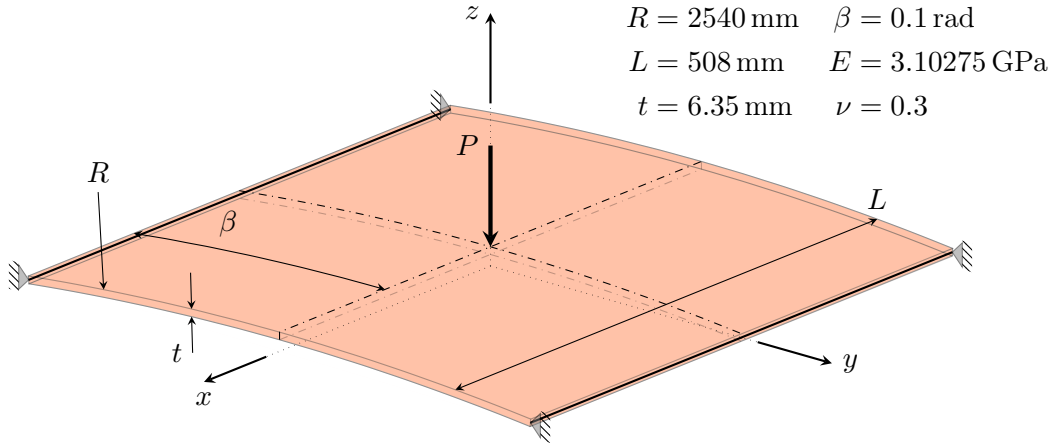

Even for linearly elastic solid materials, finite-amplitude deformation gives rise to nonlinear mechanical behavior. A well-known manifestation of such geometric nonlinearity is the presence of buckling phenomena, wherein two or more stable equilibria exist over a region of parametric hysteresis. This example considers Sabir and Lock's~\cite{sabir_lock_1972} classic benchmark of a hinged thin cylindrical section of isotropic elastic material under a central concentrated load as shown in \cref{fig:bodygeometry}. More recently, Wardle~\cite{wardle_2008} highlighted the emergence of a nontrivial symmetry-breaking buckling mode that complicates its utility as a straightforward test case for numerical solvers~\cite{sze_etal_2004,leahu-aluas_abed-meriam_2011} and suggests rich state space behavior. The present analysis will show several new bifurcations giving rise to many previously unreported solution branches, including a stable one.

Adopting a total Lagrangian approach (see Bonet, Gil \& Wood~\cite{bonet_etal_2016} for a complete derivation and further details), the nonlinear equilibrium condition for the body may be stated in the weak form as follows: find the displacement field $\bX$ in the space of virtual displacements $\mathcal{V}$ such that,
\begin{align}\label{eq:TLequation}
    \int_{\Omega_0}\boldsymbol{S}_{ij}\boldsymbol{:}\boldsymbol{\check{E}}_{ij}\,\mathrm{d}\Omega=\int_{\partial\Omega_{N,0}}\boldsymbol{f}_0\bcdot\boldsymbol{\check{X}},\qquad\forall\boldsymbol{\check{X}}\in\mathcal{V}_0,
\end{align}
where $\Omega_0\subset\mathcal{V}\in\mathbb{R}^3$ is the undeformed body, $\boldsymbol{S}_{ij}\left(\bX\right)$ are the Cartesian components of the second Piola-Kirchhoff stress tensor, $\boldsymbol{E}_{ij}\left(\bX\right)$ are the components of the Green--Lagrange strain tensor (with $\boldsymbol{\check{E}}_{ij}\left(\bX,\boldsymbol{\check{X}}\right)$ their first variation), $\partial\Omega_{N,0}$ is the Neumann part of the boundary of $\Omega$, $\boldsymbol{f}_0$ is the traction force, and $\mathcal{V}_0\subset\mathcal{V}$ is the admissible space of virtual displacements which vanish on the Dirichlet boundary $\partial\Omega_{D,0}$. Note that body forces have been neglected. This weak formulation corresponds to a strong form statement of the type in \cref{eq:equilibrium}.

Taking the point load $P$ as a parameter, the dynamics of the 3-D system are calculated using \texttt{ff-bifbox} with a FreeFEM discretization based on piecewise quadratic tetrahedral elements on a grid of $41\times41\times3$ vertices. To precisely identify symmetry breaking behavior under parameter variation, the symmetries are enforced or relaxed depending on the solution branch being traced. The symmetry group of the undeformed geometry is isomorphic to the dihedral group $D_2$. The system therefore has five possible symmetries; these are indicated as follows:
\begin{enumerate}[noitemsep]
    \item $S$, fundamental symmetry: reflections across $x=0$, $y=0$ and $180\degree$ rotations about the $z$-axis,
    \item $H_x$, half-symmetry: reflections across $x=0$,
    \item $H_y$, half-symmetry: reflections across $y=0$,
    \item $R_z$, rotational symmetry: $180\degree$ rotations about the $z$-axis, or
    \item $A$, trivial symmetry.
\end{enumerate}

\begin{figure}
    \centering
    \input{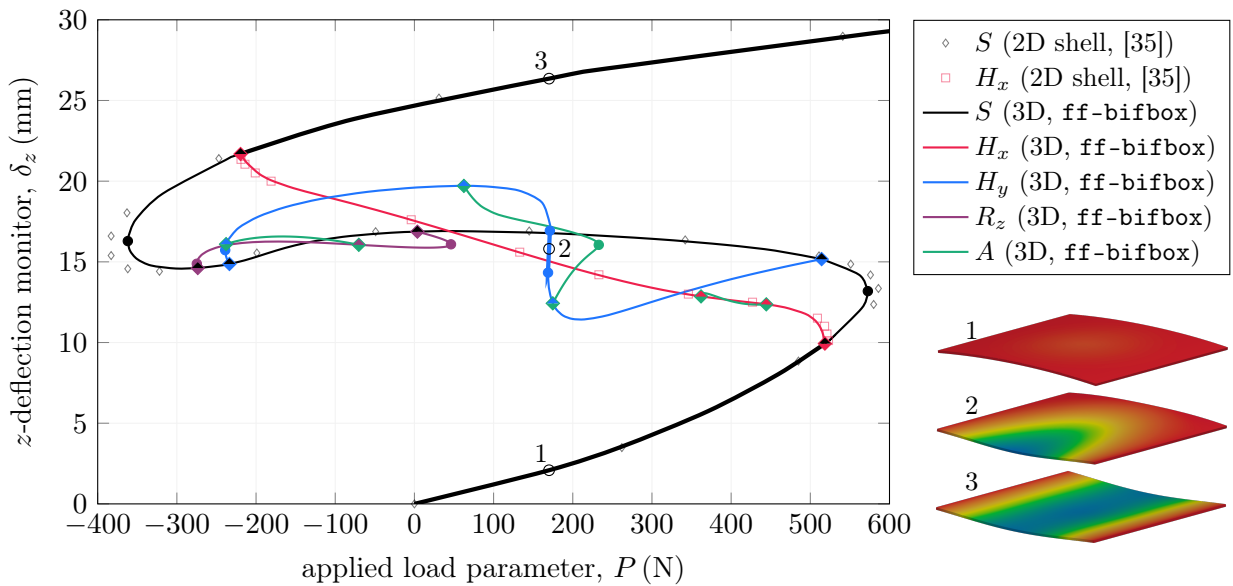}
    \caption{Bifurcation diagram showing the $z$-deflection of the center of the plate under parametric variations in $P$ and isometric visualizations of the $z$-deflection for the three stable solutions at $P=170\,N$. Linearly stable and unstable solutions are indicated by thick and thin lines, respectively. Saddle--node and symmetry breaking bifurcations are indicated on the diagram by circle and square symbols, respectively, with colors indicating the symmetry or symmetries of the bifurcating solutions. The $S$ and $H_x$ reference data from the thin shell model of Leahu-Aluas \& Abed-Meriam~\cite{leahu-aluas_abed-meriam_2011} are shown.}
    \label{fig:bucklingplate}
\end{figure}

As above, a complete example capable of reproducing the following results is available on the repository. The dynamics are summarized in the bifurcation diagram of \cref{fig:bucklingplate}. In agreement with~\cite{wardle_2008}, the results indicate a primary loss of stability as $P$ is increased that is associated with an asymmetric buckling mode. In an experimental or time-marching context, this indicates the critical loading where snap-through of the structure would spontaneously occur. However, using the branch tracing analysis, several further bifurcations associated with each possible loss of discrete symmetry of the system are also identified and traced. The extension of these broken-symmetry solutions results in a rich bifurcation diagram with multiple saddle states. Perhaps the most interesting finding from this analysis, however, is the existence of a small but finite interval of static stability with $H_y$ symmetry. This represents a stable attractor for an extremely limited range of $P$ that has not been reported in prior work. One additional quantitative observation is that the 2-D shell models used by prior works lead to a slight overprediction of the allowable loads before snap-through and snap-back in comparison to the 3-D result. A convergence test with increasing grid refinement confirmed that the result shown in \cref{fig:bucklingplate} is robust to the discretization, indicating that the discrepancy is indeed associated with the thin shell approximation. 

%%%%%%%%%%%%%%%%%%%%%%%%%%%
\subsection{Two-dimensional compressible flow past a circular cylinder}\label{ssec:cylinder}

This example demonstrates the use of \texttt{ff-bifbox} in analyzing the sequence of bifurcations leading to unsteadiness and symmetry breaking in two-dimensional sub- and transonic flow past a circular cylinder, making extensive use of FreeFEM's adaptive meshing capability. The subsonic problem has been well-studied by many authors~\cite{barkley_2006, sipp_lebedev_2007, canuto_taira_2015, fani_etal_2018, gallaire_etal_2016, mantic-lugo_etal_2014}. It has also been used more recently for validation of other open-source tools such as StabFEM~\cite{fabre_etal_2019} and BROADCAST~\cite{poulain_etal_2023}, among others. The onset of B\'{e}nard--von K\'{a}rm\'{a}n vortex shedding in the laminar flow past a cylinder is widely viewed as the prototypical example of a supercritical Hopf bifurcation in fluid mechanics. However, this analysis will show that the supercritical nature of this bifurcation is not generally preserved for compressible flow--even in the subsonic regime.

Assuming a calorically perfect gas with constant diffusivity and viscosity properties, the compressible Navier--Stokes equations governing the evolution of the velocity, temperature, and pressure fields $\bq=\left[\bu,T,p\right]^T$, may be written in dimensionless form as,
\begin{subequations}\label{eq:NSE}
\begin{align}
    \left(\gamma{}M^2p+1\right)\left(\frac{\partial\bu}{\partial t}+\bu\bcdot\bnabla\bu
    \right)+T\left(\bnabla{}p-\bnabla\bcdot\btau\right)&=0,\\
    \hspace{-1mm}\left(\gamma{}M^2p+1\right)\left(\frac{\partial{}T}{\partial t}+\bu\bcdot\bnabla{}T+\left(\gamma-1\right)T\left(\bnabla\bcdot\bu\right)\right)-\left(\gamma-1\right)\gamma{}M^2T\left(\bnabla\bu\boldsymbol{:}\btau\right)-\frac{\gamma{}T}{RePr}\nabla^2T&=0,\\
    \gamma{}M^2T\left(\frac{\partial{}p}{\partial t}+\bu\bcdot\bnabla{}p\right)-\left(\gamma{}M^2p+1\right)\left(\frac{\partial{}T}{\partial t}+\bu\bcdot\bnabla{}T-T\left(\bnabla\bcdot\bu\right)\right)&=0,
\end{align}
\end{subequations}
where the 2-D deviatoric viscous stress tensor is $\btau={Re}^{-1}\left(\bnabla\bu+\left(\bnabla\bu\right)^T-\boldsymbol{I}\left(\bnabla\bcdot\bu\right)\right)$ and the density has been eliminated using the ideal gas equation of state, $\gamma{}M^2p+1=\rho{}T$. Four dimensionless parameters---the Reynolds number $Re$, the Prandtl number $Pr$, the Mach number $M$, and the specific heat ratio $\gamma$---have been introduced based on the freestream fluid properties and the cylinder diameter. In the following, $Pr=0.72$ and $\gamma=1.4$ are fixed to focus on the influence of $Re$ and $M$. Clearly, \cref{eq:NSE} is of the form given in \cref{eq:nonlinearPDE}. Note that, with the chosen scaling and formulation of \cref{eq:NSE}, there is no singularity in the incompressible limit at $M=0$.

\begin{figure}
    \centering
    % \begin{subfigure}{0.496\textwidth}
    % \centering
    \begin{tikzpicture}[scale=0.45]
    \def\xp{12.5};
    \def\ls{2.5};
    \def\xm{5};
    \filldraw[paperorange!10] (-\xm,0) --++ (0,\xm) --++ (\xp,0) --++ (0,-\xm) -- (\xp-\xm-\ls,0) --++ (0,\xm-\ls) --++ (-\xp+2*\ls,0) --++ (0,-\xm+\ls) -- cycle;
    \filldraw[paperblue!10] (\xp-\xm-\ls,0) --++ (0,\xm-\ls) --++ (-\xp+2*\ls,0) --++ (0,-\xm+\ls) --(-0.5,0) arc (180:0:0.5) -- cycle;
    \draw[thin] (-\xm,0) --++ (0,\xm) node[midway,left=-1pt]{$\Gamma_{in}$} --++ (\xp,0) node[midway,above]{$\Gamma_{lat}$} --++ (0,-\xm) node[midway,right=-1pt]{$\Gamma_{out}$} ;
    \draw[dashed, black!75, thin] (-\xm+\ls,0) --++ (0,\xm-\ls) --++ (\xp-2*\ls,0) --++ (0,-\xm+\ls);
    \draw[dash dot] (-\xm,0) -- (-0.5,0)node[pos=0.3,below]{$\Gamma_{sym}$};
    \draw[very thick] (-0.5,0) arc (180:0:0.5)node[midway,above]{$\Gamma_{wall}$};
    \draw[dash dot] (0.5,0) -- (\xp-\xm,0)node[midway,below]{$\Gamma_{sym}$};
    \draw[-stealth] (\xp-\xm,0)--++(1,0)node[below]{$x$};
    \draw[-stealth] (0,\xm)--++(0,1)node[left]{$y$};
    \draw[dotted] (0,\xm)--++(0,-3.5);
    \draw[black!75,stealth-stealth] (\xm-\ls,\xm)--++(0,-\ls)node[midway,right]{$l_s$};
    \filldraw (-0.5,0) circle (0.1) node[below left=-1pt and -9pt]{\tiny{$(-\frac{1}{2},0)$}};
    \filldraw (0.5,0) circle (0.1) node[below right=-1pt and -9pt]{\tiny{$(\frac{1}{2},0)$}};
    \filldraw (-\xm,\xm) circle (0.1) node[above]{\tiny{$(-100,100)$}};
    \filldraw (\xp-\xm,\xm) circle (0.1) node[above]{\tiny{$(150,100)$}};
    \filldraw[black!75] (-\xm+\ls,\xm-\ls) circle (0.1) node[above]{\tiny{$(-20,20)$}};
    \filldraw[black!75] (\xp-\xm-\ls,\xm-\ls) circle (0.1) node[above]{\tiny{$(50,20)$}};
    \end{tikzpicture}
    \caption{Illustration (not to scale) of the flow configuration over the circular cylinder. The blue interior region corresponds to the physical domain of interest and the orange outer region denotes the absorbing layer.}
    \label{fig:flowconfig}
\end{figure}
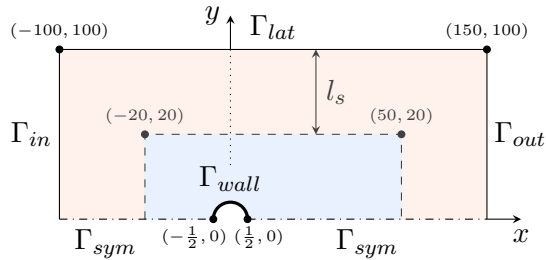

A schematic of the flow configuration is shown in \cref{fig:flowconfig}. With $\bn$ denoting the outward normal unit vector on the boundary, the flow is configured with adiabatic no-slip walls ($\Gamma_{wall}$: $\bu=0$, $\bn\bcdot\bnabla{}T=0$), uniform isothermal inflow conditions ($\Gamma_{in}$: $\left[u_x,u_y,T\right]^T=\left[1,0,1\right]^T$), free-slip adiabatic lateral boundaries ($\Gamma_{lat}$: $\bn\bcdot\bu=0$, $\bn\times\left(\btau\bcdot\bn\right)=0$, $\bn\bcdot\bnabla{}T=0$), and stress-free adiabatic outflow conditions ($\Gamma_{out}$: $\btau\bcdot\bn-p\bn=0$, $\bn\bcdot\bnabla{}T=0$). Symmetry conditions are used along the flow centerline $\Gamma_{sym}$ to reduce computational effort and allow precise identification of transverse symmetry breaking. Note that, as in prior compressible flow studies~\cite{rowley_etal_2002,colonius_2004,yamouni_etal_2013,fani_etal_2018}, the physical part of the computational domain is surrounded by a thick absorbing layer to mitigate boundary scattering effects. Within the absorbing layer, artificial damping terms given by $\beta{}_sd_s^2\left(\bx\right)\left[u_x-1, u_y, \gamma{}\left(T-1\right),\gamma{}M^2p\right]^T$ 
are added to the residual of \cref{eq:NSE}, where $\beta_s=10^{-4}$ is a small constant and $d_s\left(\bx\right)$ is the Euclidean distance to the physical region. The in-plane discretization of the computational domain consists of an adaptively-refined triangulation and projection onto a Taylor--Hood-type mixed space consisting of piecewise quadratic (i.e. $\mathbb{P}_2$) elements for $\bu$ and $T$ and piecewise linear (i.e. $\mathbb{P}_1$) elements for $p$. Adaptive refinement is based on a $\mathbb{P}_1$ interpolation error of 1\% for the global solution using FreeFEM's \texttt{adaptmesh} function~\cite{hecht_2012}. For the considered range of parameters, this results in meshes with $O\left(10^4\right)$ vertices and states with $O\left(10^5\right)$ discrete degrees of freedom.

\begin{figure}[t!]
    \centering
    \input{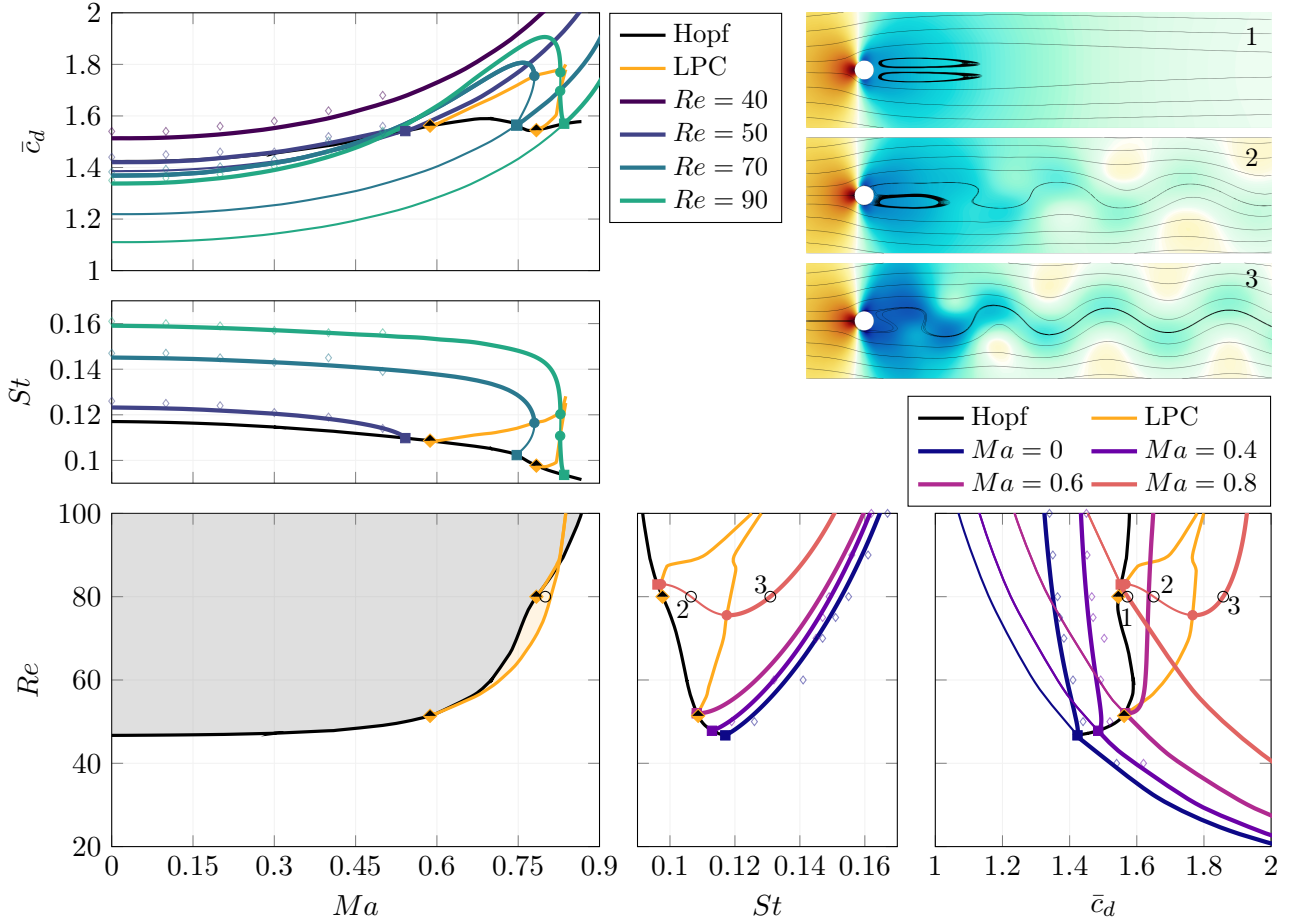}
    \caption{Bifurcation diagrams of $\bar{c}_d$ and $St$ versus $M$ and $Re$, $M$--$Re$ stability map, and flow visualization summarizing the dynamics of the 2-D compressible flow past a circular cylinder. The projections of the Hopf bifurcation curve associated with time-periodic vortex shedding are plotted in black, with the gray shaded portion of the stability map indicating linear instability of the steady state. The two-tone diamonds indicate codimension-2 Bautin bifurcation points. The projections of the LPC curves are drawn in yellow, with the yellow shaded area in the stability map indicating the bistable region. In the bifurcation diagrams, Hopf and LPC points are indicated with square and round markers, respectively, with asymptotically stable and unstable solutions indicated by thick and thin lines, respectively. Bifurcation diagrams are shown for $Re=[40,50,70,90]$ over varying $M$ and $M=[0,0.4,0.6,0.8]$ over varying $Re$. Reference data from direct numerical simulations of Canuto \& Taira~\cite{canuto_taira_2015} are shown by color-matched diamond markers. Snapshots of the stable steady state, saddle periodic solution, and stable periodic solution at $(M,Re)=(0.8,80)$ are visualized via instantaneous pressure contours and velocity streamlines over $(x,y)\in[-3,-3]\times[21,3]$.}
    \label{fig:bifmaps}
\end{figure}

Bifurcation diagrams of the time-mean drag coefficient $\bar{c}_d$ and Strouhal number $St$, a $Re$--$M$ stability map, and flow visualizations at $(M,Re)=(0.8,80)$ are given in \cref{fig:bifmaps} to represent the system's nonlinear dynamics. In agreement with prior work~\cite{canuto_taira_2015,fani_etal_2018}, the flow exhibits a steady symmetric solution for sufficiently low $Re$ over all considered $M$. With increasing $Re$, the system undergoes a Hopf bifurcation along the neutral curve shown in \cref{fig:bifmaps}, signifying the intersection of a limit cycle solution manifold with the steady symmetric manifold. On the limit cycle, the flow exhibits a time-periodic, two-dimensional B\'enard--von K\'arm\'an state which breaks the bilateral symmetry. As in prior work~\cite{sipp_lebedev_2007,mantic-lugo_etal_2014,fani_etal_2018}, the bifurcation to this limit cycle is supercritical at low $M$. However, the results also indicate a previously undescribed transition to subcritical behavior mediated by a codimension-2 Bautin bifurcation at $M\approx0.59$, which is detected automatically by \texttt{ff-bifbox} during continuation of the Hopf curve. A second Bautin bifurcation occurs further along the neutral curve, indicating a more complex ``pleated'' structure of the periodic solution branch with bistable limit cycle states. Beyond the second Bautin bifurcation, the extent of the bistable regime is delimited by the loci of limit points of cycles (LPC) emanating from the Bautin points. Visualizations of the bistable steady and time-periodic states, as well as the unstable time-periodic saddle state, at $(M,Re)=(0.8,80)$ are shown in \cref{fig:bifmaps}. Overall, the results indicate far richer nonlinear dynamics in the transonic regime than has been revealed by prior investigations at lower $M$.

%%%%%%%%%%%%%%%%%%%%%%%%%%%
\section{Conclusion} \label{sec:conclusion}
%%%%%%%%%%%%%%%%%%%%%%%%%%%
This article introduces \texttt{ff-bifbox}, a new open-source branch tracing and bifurcation analysis toolbox specialized for large, sparse, nonlinear PDEs and systems of nonlinear PDEs. The toolbox takes advantage of finite element discretization and mesh adaptation in FreeFEM~\cite{hecht_2012} and efficient, scalable numerical linear algebra routines in PETSc~\cite{balay_etal_2026} and SLEPc~\cite{hernandez_etal_2005}. Here, details of the underlying theory and approach are presented, along with example implementations and results for three common classes of nonlinear PDEs: a reaction--diffusion system, a solid nonlinear elasticity system, and a compressible fluid flow system. These results are used for validation of the implementation, but also describe novel findings for each system.

The \texttt{ff-bifbox} code is under ongoing development to extend its capabilities across a broader scope of dynamical systems analysis, to improve parallel performance, and to make it more robust and flexible for complex multiphysics problems. Some noteworthy areas of primary focus include: improved capabilities for detection, computation, and continuation of local bifurcations of limit cycle solutions; direct computation and normal form evaluation for codimension-2 Bautin, cusp, and Bogdanov--Takens bifurcations; and harmonic resolvent analysis for periodic orbits. Further development efforts toward robust implementations of deflation-based iterative techniques for near-singular systems near bifurcation points and efficient preconditioners for analysis of periodic orbits based on the harmonic balance method~\cite{moulin_2020, sierra_etal_2021, sierra-ausin_etal_2022} are also underway.

\subsection*{CRediT authorship contribution statement}
\textbf{C.~Douglas:} Conceptualization, Data curation, Formal analysis, Investigation, Methodology, Resources, Software, Validation, Visualization, Writing -- original draft, Writing -- review and editing. \textbf{P.~Jolivet:} Methodology, Software, Validation, Writing -- review and editing.

\subsection*{Declaration of competing interest}
The authors declare that they have no known competing financial interests or personal relationships that could have appeared to influence the work reported in this paper.

\subsection*{Data availability}
The data used in this manuscript may be freely reproduced using the complete source code available in the supplementary materials and at \href{https://github.com/cmdoug/ff-bifbox}{https://github.com/cmdoug/ff-bifbox} under the GPL-3.0 license. The authors have also applied a CC-BY public copyright license to the present manuscript and all subsequent revisions up to the Author Accepted Manuscript.

\printbibliography
\end{document}